\documentclass[12pt]{amsart}

\usepackage{graphicx,amssymb,latexsym,amsfonts,txfonts,amsmath,amsthm}
\usepackage{pdfsync,color,tabularx,rotating}
\usepackage[all,cmtip]{xy}

\usepackage{tikz}

\advance\hoffset by -0.9 truecm

\newtheorem{theo}{Theorem}

\newtheorem{cor}[theo]{Corollary}
\newtheorem{lemma}[theo]{Lemma}

\theoremstyle{remark}

\theoremstyle{remark}

\theoremstyle{remark}

\begin{document}

\title{Groups of automorphisms of Riemann surfaces and maps of genus $p+1$ where $p$ is prime}

\author{Milagros Izquierdo, Gareth A. Jones and Sebasti\'an Reyes-Carocca}

\address{Department of Mathematics, Link\"oping University, 58183 Link\"oping, Sweden}
\email{milagros.izquierdo@liu.se}

\address{School of Mathematical Sciences, University of Southampton, Southampton SO17 1BJ, UK}
\email{G.A.Jones@maths.soton.ac.uk}

\address{Departamento de Matem\'atica y Estad\'\i stica, Universidad de La Frontera, Avenida Francisco Salazar 01145, Temuco, Chile}
\email{sebastian.reyes@ufrontera.cl}



\maketitle


\begin{abstract}
We classify compact Riemann surfaces of genus $g$, where $g-1$ is a prime $p$, which have a group of automorphisms of order $\rho(g-1)$ for some integer $\rho\ge 1$, and determine isogeny decompositions of the corresponding Jacobian varieties. This extends results of Belolipetzky and the second author for $\rho>6$, and of the first and third authors for $\rho=3, 4, 5$ and $6$. As a corollary we classify the orientably regular hypermaps (including maps) of genus $p+1$, together with the non-orientable regular hypermaps of characteristic $-p$, with automorphism group of order divisible by the prime $p$; this extends results of Conder, \v Sir\'a\v n and Tucker for maps.
\end{abstract}


\section{Introduction}\label{intro}

A compact Riemann surface $\mathcal S$ of genus $g\ge 2$ has a finite automorphism group, of order at most $84(g-1)$. It is well known that for a given genus $g$ the possibilities for surfaces $\mathcal S$ and their automorphism groups depend heavily on the factorisation of the Euler characteristic $\chi=2-2g$, since divisors of $\chi$ allow such surfaces to occur as unbranched coverings of those of smaller genus. From this point of view, the simplest case to consider is therefore that in which $g-1$ is a prime~$p$. In~\cite{BJ}, Belolipetzky and the second author considered this situation on the assumption that $\mathcal S$ has a group $G$ of automorphisms of order $\rho(g-1)$ where $\rho\ge \lambda$ for some $\lambda>6$; they showed that if $p$ is sufficiently large as a function of $\lambda$ (to avoid finitely many sporadic cases) then $\mathcal S$ and $G$ lie in one of six infinite families, each with a simple construction. This work has been reinterpreted and taken further in the context of orientably regular maps by Conder, \v Sir\'a\v n and Tucker in~\cite{CST}. More generally, Conder and Kulkarni~\cite{CK} have investigated sequences of groups of automorphisms of order $ag+b$ for constants $a, b\in{\mathbb Q}$, such as the Accola--Maclachlan groups of orders $8(g+1)$ and $8(g+3)$.

The group-theoretic techniques available to study this problem divide it naturally into two general cases, according to whether or not $p$ divides $|G|$ (or equivalently $\rho\in{\mathbb Z}$). If it does, then the Sylow and  Schur--Zassenhaus theorems imply that, for all but finitely many primes $p$, $G$ is a semidirect product $P\rtimes Q$ of a normal Sylow $p$-subgroup $P\cong C_p$ by a group $Q$ of order~$\rho$. Moreover, $Q\le{\rm Aut}\,{\mathcal T}$ where ${\mathcal T}={\mathcal S}/P$ is a Riemann surface of genus $2$, so the possibilities for $Q$ are limited. The small number of exceptional primes $p$ can be dealt with by ad hoc methods. (The primes $p=2, 3$ and $5$ are often omitted, since automorphism groups for very small genera $g$ behave differently and are well-known: see~\cite{Bro, Con}, for example.) The results in~\cite{BJ} for $\rho>6$ have recently been extended by the first and third authors~\cite{IR, Rey} to the cases where $\rho=3, 4, 5$ or $6$ and $p\ge 7$; for instance, they show that for $\rho=5$ there are no new surfaces (the groups $G$ with $\rho=5$ are all subgroups of those appearing in~\cite{BJ} as full automorphism groups of surfaces $\mathcal S$ with $\rho=10$), whereas for $\rho=6$ they describe an infinite family of surfaces, two members of which also appear in~\cite{BJ} with larger automorphism groups; this family also realises the groups arising for $\rho=3$.

Here we will build on these results, filling in gaps (such as for small primes where $\rho>6$) to give a unified treatment and a complete classification of the groups $G$ and surfaces $\mathcal S$ for all integers $\rho\ge 1$. We also describe some group actions, such as those in case~(v) of Theorem~\ref{integerthm}(a), which are only implicit in~\cite{IR}, where the emphasis is more on the surfaces than on the groups. For conciseness of exposition and proof, the main result, Theorem~\ref{integerthm}, is stated (below) only for integers $\rho\ge 3$ and primes $p\ge 7$. Small values of $\rho$ and $p$, which lead to less uniform behaviour, are discussed separately towards the end of the paper.

In order to state our results, we introduce some notation. For each prime $p>2$ and divisor $r$ of $p-1$ let us define a group
\[G_{p,r}:=\langle a, b\mid a^p=b^r=1,\, bab^{-1}=a^{\omega}\rangle\]
where $\omega$ is a primitive $r$th root of $1$ in $\mathbb Z_p$. This is a semidirect product of $\langle a\rangle\cong C_p$ by $\langle b\rangle\cong C_r$, with the latter acting faithfully by conjugation on the former. Up to isomorphism this group, denoted by $C_p\rtimes_rC_r$ in~\cite{IR}, is independent of the choice of $\omega$, and is the unique subgroup of order $pr$ in the affine group $AGL_1(p)\cong G_{p,p-1}$. For example, $G_{p,2}$ is a dihedral group $D_p$ of order $2p$.

Each compact Riemann surface $\mathcal S$ of genus $g\ge 2$ is isomorphic to a quotient ${\mathbb H}/K$ of the hyperbolic plane $\mathbb H$ by a surface group $K$ of genus $g$. A group $G$ (necessarily finite) is isomorphic to a subgroup of ${\rm Aut}\,\mathcal S$ if and only if $G\cong\Gamma/K$ where $\Gamma$ is a cocompact Fuchsian group containing $K$ as a normal subgroup. The signature of $\Gamma$ has the form
\[\sigma=(\gamma;m_1,\ldots, m_k)\]
for some genus $\gamma$ (of ${\mathbb H}/\Gamma\cong {\mathcal S}/G$) and elliptic periods $m_i\ge 2$.
We will also refer to $\sigma$ as the signature of the action of $G$ on $\mathcal S$. The order of the periods $m_i$ is irrelevant. For brevity, we will omit $\gamma$ and write simply $(m_1,\ldots, m_k)$ in the (rather frequent) cases where $\gamma=0$, and we will denote an elliptic period $m$ repeated $r$ times by $m^{[r]}$. Our main result is the following:

\begin{theo}\label{integerthm}
Let $\mathcal S$ be a compact Riemann surface of genus $g=p+1$ for some prime $p\ge 7$. 

\smallskip

\noindent {\rm (a)} There is a subgroup $G\le{\rm Aut}\,{\mathcal S}$ of order $|G|=\rho(g-1)$ for some integer $\rho\ge 3$ if and only if one of the following holds, where $\sigma$ denotes the signature of the action of $G$:
\begin{enumerate}
\item[(i)] $\rho=12$, $\sigma=(2,6,6)$ and $G\cong G_{p,6}\times C_2$ where $p\equiv 1$ {\rm mod}~$(3)$;
\item[(ii)] $\rho=10$, $\sigma=(2,5,10)$ and $G\cong G_{p,10}$ where $p\equiv 1$ {\rm mod}~$(5)$;
\item[(iii)] $\rho=8$, $\sigma=(2,8,8)$ and $G\cong G_{p,8}$ where $p\equiv 1$ {\rm mod}~$(8)$;
\item[(iv)] $\rho=6$, $\sigma=(3,6,6)$ and $G\cong G_{p,6}$ or $G_{p,3}\times C_2$ where $p\equiv 1$ {\rm mod}~$(3)$;
\item[(v)] $\rho=6$, $\sigma=(2,2,3,3)$ and $G\cong G_{p,6}$ where $p\equiv 1$ {\rm mod}~$(3)$;
\item[(vi)] $\rho=5$, $\sigma=(5,5,5)$ and $G\cong G_{p,5}$ where $p\equiv 1$ {\rm mod}~$(5)$;
\item[(vii)]  $\rho=4$, $\sigma=(2,2,4,4)$ and $G\cong G_{p,4}$ where $p\equiv 1$ {\rm mod}~$(4)$;
\item[(viii)] $\rho=4$, $\sigma=(2^{[5]})$ and $G\cong G_{p,2}\times C_2\cong D_{2p}$;
\item[(ix)] $\rho=3$,  $\sigma=(3^{[4]})$ and $G\cong G_{p,3}$ where $p\equiv 1$ {\rm mod}~$(3)$;
\item[(x)] $\rho=84$, $\sigma=(2,3,7)$ and $G\cong PSL_2(13)$ where $p=13$;
\item[(xi)] $\rho=48$, $\sigma=(2,3,8)$ and $G\cong PGL_2(7)$ where $p=7$;
\item[(xii)] $\rho=24$, $\sigma=(3,3,4)$ and $G\cong PSL_2(7)$ where $p=7$.
\end{enumerate}

\smallskip

\noindent{\rm (b)} In cases (i) and (iv), for each $p$ there is a single chiral pair of surfaces $\mathcal S_1$ and $\overline{\mathcal S}_1$, the same in each case, each surface admitting both groups $G$ in (iv). In cases~(ii) and (vi) there are two chiral pairs $\mathcal S_2, \overline{\mathcal S}_2$ and $\mathcal S_2', \overline{\mathcal S'_2}$, the same in each case. In case (iii) there is a single chiral pair of surfaces $\mathcal S_3, \overline{\mathcal S}_3$. In case (v) there is an infinite family of surfaces, of real dimension $2$, with an action of $G_{p,6}$ of signature $(2,2,3,3)$, restricting to the action of $G_{p,3}$ in (ix); two of these surfaces, namely the chiral pair $\mathcal S_1, \overline{\mathcal S}_1$ in (i), also both admit actions of $G_{p,6}$ and $G_{p,3}\times C_2$ of signature $(3,6,6)$ in (iv). In cases (vii) and (viii) there is an infinite family of surfaces for each $p$, of real dimension $2$ and $4$ respectively; if $p\equiv 1$ {\rm mod}~$(8)$ then two of these in (vii) are the chiral pair $\mathcal S_3, \overline{\mathcal S}_3$ in (iii). In case~(x) there are three surfaces, and in cases~(xi) and (xii) there are two, the same in each case. 

\smallskip

\noindent{\rm (c)} The surfaces $\mathcal S$ in case~(iv) and a chiral pair of those in case~(v) are the surfaces $\mathcal S_1, \overline{\mathcal S}_1$ in (i) with automorphism group $A:={\rm Aut}\,{\mathcal S}\cong G_{p,6}\times C_2$; when $p=7$ two of those in case~(v) have $A\cong PGL_2(7)$ in (xi), and when $p=13$ three of those in case~(v) have $A\cong PSL_2(13)$ in (x). The surfaces in case~(vi) have $A\cong G_{p,10}$ in (ii). In case~(vii), if $p\equiv 1$ {\rm mod}~$(8)$ a chiral pair $\mathcal S_3, \overline{\mathcal S}_3$ have $A\cong G_{p,8}$ in (iii). In case~(viii), if $p\equiv 1$ {\rm mod}~$(3)$ a chiral pair $\mathcal S_1, \overline{\mathcal S}_1$ have $A\cong G_{p,6}\times C_2$ in (i). The surfaces in case~(ix) are the infinite family in case~(v), with $G_{p,3}$ acting as a subgroup of index $2$ in $G_{p,3}\times C_2$, and with automorphism groups as described here for case~(v). The surfaces in case~(xii) have $A\cong PGL_2(7)$ in (xi). All other surfaces $\mathcal S$ have automorphism group $A=G$. 
\end{theo}

This theorem confirms and extends results obtained earlier by Belolipetzky and the second author~\cite{BJ} for $\rho>6$, and more recently by the first and third authors~\cite{IR, Rey} for $\rho=3,4,5,6$. There are similar but less uniform results for primes $p\le 5$ and for $\rho=1$ and $2$, discussed briefly in Sections~\ref{smallp} and \ref{smallrho} after the proof of Theorem~\ref{integerthm}.

The Jacobian variety $J\mathcal S$ of a compact Riemann surface $\mathcal S$ of genus $g\ge 2$ is a principally polarized abelian variety of dimension $g$, namely a complex torus which is projective.  It is known that Jacobians are irreducible, in the sense that they are not isomorphic to the product of two abelian subvarieties of lower dimension. The relevance of Jacobian varieties lies in part in Torelli's theorem, that two compact Riemann surfaces are isomorphic if and only if their Jacobians are isomorphic as principally polarized abelian varieties.

Let $G$ be a finite group acting conformally on a compact Riemann surface $\mathcal S$. It is well known that this action induces an action of $G$ on $J\mathcal S$ and this, in turn, induces an isogeny decomposition which is $G$-equivariant (see~\cite{CR, LR}). This decomposition of Jacobians under group actions has been extensively studied, following early papers by Wirtinger~\cite{Wir} and Schottky and Jung~\cite{SJ}. For decomposition of Jacobians with respect to specific groups, see~\cite{CRR, FLP, HJQR, IJR, PR, RR98, R19, RR20}.

The following result extends and confirms previous results in~\cite{IR} and~\cite{Rey} and provides a complete treatment of isogeny decompositions for each surface $\mathcal S$ in Theorem~\ref{integerthm}.

\begin{theo}\label{Jacobians}
For each surface $\mathcal S$ in Theorem~\ref{integerthm}, the Jacobian $J\mathcal S$ decomposes up to isogeny as follows. For $\mathcal S$ in case~(i) we have
\[J{\mathcal S}\sim J{\mathcal T}\times(J{\mathcal C})^6\]
where ${\mathcal T}={\mathcal S}/\langle a\rangle$ has genus $2$ and ${\mathcal C}={\mathcal S}/\langle bt\rangle$ has genus $(p-1)/6$ with $t$ generating the direct factor $C_2$ of $G_{p,6}\times C_2$.
In case~(iv) with an action of $G_{p,3}\times C_2$ we have
\[J{\mathcal S}\sim J{\mathcal T}\times(J{\mathcal C})^3\]
where ${\mathcal T}={\mathcal S}/\langle a\rangle$ has genus $2$ and ${\mathcal C}={\mathcal S}/\langle b\rangle$ has genus $(p-1)/3$.
In case~(viii) we have
\[J{\mathcal S}\sim E\times J{\mathcal C}_1\times J{\mathcal C}_2\]
where $E={\mathcal S}/\langle at\rangle$ is an elliptic curve, ${\mathcal C}_1={\mathcal S}/\langle bt\rangle$ has genus $(p-1)/2$ and ${\mathcal C}_2={\mathcal S}/\langle ab\rangle$ has genus $(p+1)/2$.
In the remaining cases, with an action of $G_{p,r}$, we have
\[J{\mathcal S}\sim J{\mathcal T}\times(J{\mathcal C})^r\]
where ${\mathcal T}={\mathcal S}/\langle a\rangle$ has genus $2$ and ${\mathcal C}={\mathcal S}/\langle b\rangle$ has genus $(p-1)/r$.
\end{theo}

Here case~(iii) is new, cases~(iv) and (ix) are new but implicit in~\cite[Theorem 2]{IR} since they are contained in cases~(i) and (v) respectively, case~(vi) is new but implicit in~\cite[Theorem 1]{IR} since it is contained in case~(ii), and cases~(vii) and (viii) are dealt with in~\cite[Theorem 3]{Rey}. The proof of Theorem~\ref{Jacobians} is given in Section~\ref{JacProof}.

Connections of these results with maps and hypermaps are discussed in Section~\ref{maps} where, as a corollary to Theorem~\ref{integerthm}, we obtain the following classification, where the numbering of cases follows and refers to that in Theorem~\ref{integerthm}(a):

\begin{theo}\label{mapsthm}
The orientably regular maps or hypermaps of genus $g=p+1$ for some prime $p\ge 7$, with orientation-preserving automorphism group $G$ of order divisible by $p$, are as follows (up to duality or triality, permuting the roles of vertices, edges and faces):
\begin{itemize}
\item[{\rm(i)}] for $p\equiv 1$ {\rm mod}~$(3)$ the surfaces ${\mathcal S}_1$ and $\overline{\mathcal S}_1$ in Theorem~\ref{integerthm}(a)(i) support a chiral pair of orientably regular maps of type $\{6,6\}$ with $G\cong G_{p,6}\times C_2$;
\item[{\rm(ii)}] for $p\equiv 1$ {\rm mod}~$(5)$ the surfaces ${\mathcal S}_2$, $\overline{\mathcal S}_2$, ${\mathcal S}'_2$ and $\overline{\mathcal S'}_2$ in Theorem~\ref{integerthm}(a)(ii) support two chiral pairs of orientably regular maps of type $\{5,10\}$ with $G\cong G_{p,10}$;
\item[{\rm(iii)}] for $p\equiv 1$ {\rm mod}~$(8)$ the surfaces ${\mathcal S}_3$ and $\overline{\mathcal S}_3$ in Theorem~\ref{integerthm}(a)(iii) support a chiral pair of orientably regular maps of type $\{8,8\}$ with $G\cong G_{p,8}$;
\item[{\rm(iv)}] for $p\equiv 1$ {\rm mod}~$(3)$ the surfaces ${\mathcal S}_1$ and $\overline{\mathcal S}_1$ in Theorem~\ref{integerthm}(a)(i) and (iv) support two chiral pairs of orientably regular hypermaps of type $(3,6,6)$, one each with $G\cong G_{p,6}$ or $G_{p,3}\times C_2$;
\item[{\rm(vi)}] for $p\equiv 1$ {\rm mod}~$(5)$ the surfaces ${\mathcal S}_2$, $\overline{\mathcal S}_2$, ${\mathcal S}'_2$ and $\overline{\mathcal S'}_2$ in Theorem~\ref{integerthm}(a)(ii) and (vi) support twelve orientably regular hypermaps of type $(5,5,5)$ with $G\cong G_{p,5}$;
\item[{\rm(x)}] for $p=13$ the three surfaces $\mathcal S$ in Theorem~\ref{integerthm}(a)(x) support three regular maps of type $\{3,7\}$ with $G\cong PSL_2(13)$, one on each surface; 
\item[{\rm(xi)}] for $p=7$ the two surfaces $\mathcal S$ in Theorem~\ref{integerthm}(a)(xi) support two regular maps of type $\{3,8\}$ with $G\cong PGL_2(7)$, one on each surface;
\item[{\rm(xii)}] for $p=7$ the two surfaces $\mathcal S$ in Theorem~\ref{integerthm}(a)(xi) support two regular hypermaps of type $(3,3,4)$ with $G\cong PSL_2(7)$, one on each surface. 
\end{itemize}
\end{theo}

As before, primes $p\le 5$ are omitted for conciseness, but are easily dealt with. This theorem can also be regarded as a classification of the regular dessins d'enfants (see~\cite{JW}) satisfying the same conditions on their genus and automorphism group. There is a similar classification of non-orientable regular maps and hypermaps of characteristic $-p$ in Section~\ref{nonor}. These results extend to hypermaps some earlier results for maps by Conder, \v Sir\'a\v n and Tucker in~\cite{CST}, where they also consider the case where $p$ does not divide $|G|$. For small $p$, these maps and hypermaps are identified in Sections~\ref{maps} and \ref{nonor} with the corresponding entries in Conder's computer-generated lists~\cite{Con}.

In Theorem~\ref{integerthm} there is an obvious contrast between cases~(i) to (ix), which describe infinite sequences (guaranteed by Dirichlet's Theorem on primes in arithmetic progressions) exhibiting uniform behaviour, and cases~(x) to (xii) where we have small sporadic examples exhibiting irregular behaviour.  In Sections~\ref{maps} and \ref{nonor} we see the same contrast concerning the maps and hypermaps on these surfaces. This distinction, a common phenomenon for both finite groups and compact Riemann surfaces in general, is explained here by the fact (see the second paragraph above) that in cases~(i) to (ix) the Sylow $p$-subgroup $P\cong C_p$ of $G$ is normal, implying that $\mathcal S$ is a regular unbranched $p$-sheeted covering of a Riemann surface ${\mathcal T}={\mathcal S}/P$ of genus $2$, whereas in cases~(x) to (xii) $P$ is not normal in $G$, and $\mathcal S$ does not have this structure.

It is worth emphasising that this paper does not consider values $\rho\in{\mathbb Q}\setminus{\mathbb Z}$, where the methods and results (see~\cite{BJ, CST}, for instance) are different: for example, the Accola--Maclachlan groups of order $8(g+3)$ and $8(g+1)$ and their associated surfaces play an important role there, and only the first of these, with $g=3$, has $g-1$ prime and $\rho\in{\mathbb Z}$. It is hoped to revisit this situation later.

\medskip

\noindent{\bf Acknowledgment.} The authors are grateful to David Singerman for a number of helpful remarks concerning Fuchsian groups. The third author was partially supported by Fondecyt Grants 11180024 and 1190991. The first and third authors were partially supported by Redes Grant 170071.


\section{Preliminaries}\label{prelim}

Throughout this paper we will be concerned with cocompact Fuchsian groups. A {\sl Fuchsian group\/} is a discrete subgroup $\Gamma$ of $PSL_2({\mathbb R})$, acting discontinuously by M\"obius transformations on the upper half plane $\mathbb H$. We say that $\Gamma$ is {\sl cocompact\/} if the quotient space ${\mathbb H}/\Gamma$ is compact, in which case it is known that $\Gamma$ is finitely generated and contains no parabolic elements: the non-identity elements of $\Gamma$ are all elliptic and of finite order, with one fixed point in $\mathbb H$, or hyperbolic and of infinite order, with no fixed points in $\mathbb H$ and two on its boundary ${\mathbb P}^1(\mathbb R)={\mathbb R}\cup\{\infty\}$. For background on compact Riemann surfaces and Fuchsian groups, we refer to~\cite{FK}.

For the rest of this paper, $\mathcal S$ will denote a compact Riemann surface of genus $g\ge 2$. By the Uniformization Theorem, $\mathcal S$ is conformally equivalent (isomorphic) to the quotient $\mathbb{H}/K$ of $\mathbb H$ by a Fuchsian group $K$ isomorphic to the fundamental group $\Pi_g$ of $\mathcal S$. 


\subsection{Group actions, topological equivalence}

We say that a group $G$ acts on $\mathcal S$ if there is a group monomorphism $\psi: G\to {\rm Aut}\,{\mathcal S}$. The condition $g\ge 2$ implies that ${\rm Aut}\,\mathcal S$ is finite, so $G$ acts discontinuously on $\mathcal S$. The orbit-space ${\mathcal S}/G$ induced by this action $\psi$ then has a natural Riemann surface structure so that the projection ${\mathcal S} \to {\mathcal S}/G$ is holomorphic.

A group $G$ acts on $\mathcal S$ if and only if there is a Fuchsian group $\Gamma$ containing $K$ with an epimorphism $\theta: \Gamma \to G$ such that $\ker\theta=K$ (see~\cite{Sin70}); in this case ${\mathcal S}/G\cong \mathbb{H}/\Gamma$. We call $\theta$ a {\em surface epimorphism}. Since $\mathcal S$ is compact, so is ${\mathcal S}/G$, so $\Gamma$ is cocompact. It then follows that $\Gamma$ has a presentation with generators
\[A_j, B_j\quad(j=1,\ldots, \gamma)\quad{\rm and}\quad X_i\quad(i=1,\ldots, k)\]
(respectively hyperbolic and elliptic), and defining relations
\[\prod_{j=1}^{\gamma}[A_j,B_j].\prod_{i=1}^kX_i=X_i^{m_i}=1.\]
(We will use this notation for generators throughout the paper.) Here $\gamma$ is the genus of ${\mathcal S}/G$, and the elliptic periods $m_i$ indicate the order of branching at the branch points of the covering ${\mathcal S}\to{\mathcal S}/G$. The order of the elliptic periods is irrelevant, and we will usually assume that $2\le m_1\le \cdots\le m_k$.

One can encode this presentation by saying that $\Gamma$ has {\em signature}
\[\sigma=(\gamma;m_1,\ldots, m_k).\]
We will also refer to $\sigma$ as the signature of the action $\psi$ of $G$ on $\mathcal S$ (with monodromy $\theta$), or more concisely (but less precisely) as the signature of $G$. More generally, we will write $\Gamma(\sigma)$ to denote any Fuchsian group with this signature.

Two actions $\psi_1, \psi_2: G \to {\rm Aut}\,\mathcal S$ of $G$ on $\mathcal S$ are said to be {\it topologically equivalent} if there exist an automorphism $\omega$ of $G$ and an orientation-preserving self-homeomorphism $h$ of $\mathcal S$ such that
\begin{equation}\label{equivalentactions}
\psi_2(g) = h \psi_1(\omega(g)) h^{-1} \quad \mbox{for all} \; g\in G.
\end{equation}
In this case $\psi_1$ and $\psi_2$ have the same signature. Each orientation-preserving homeomorphism $h$ satisfying (\ref{equivalentactions}) yields a group automorphism $h^*$ of $\Gamma$. We let $\mathfrak{B}$ denote the subgroup of ${\rm Out}\,(\Gamma)$ consisting of the images of such automorphisms $h^*$. Equivalently, surface epimorphisms $\theta_1, \theta_2 : \Gamma \to G$ define topologically equivalent actions if and only if $\theta_2 = \omega\circ\theta_1 \circ h^*$ for some $\omega \in {\rm Aut}(G)$ and $h^* \in \mathfrak{B}$ (see~\cite{Bro90, Har71, MS}). If ${\mathcal S}/G$ has genus $\gamma=0$ then $\mathfrak B$ is generated by the braid transformations.


\subsection{Equisymmetric stratification}

Each Riemann surface $\mathcal S$ of genus $g\ge 2$ is uniformised by a surface Fuchsian subgroup $K\cong \Pi_g$ of $PSL_2({\mathbb R})$. Two subgroups uniformise isomorphic surfaces if and only if they are conjugate in $PSL_2({\mathbb R})$. We define the {\em Teichm\"uller space} ${\mathbb T}_g$ to be the quotient of the space of such embeddings $r:\Pi_g\to K\le PSL_2({\mathbb R})$ modulo conjugation in $PSL_2({\mathbb R})$; it is homeomorphic to a ball of dimension $6g-6$. The {\em modular group\/} or mapping class group ${\rm Mod}_g:={\rm Aut}^+(\Pi_g)/{\rm Inn}(\Pi_g)$ acts by composition on $\mathbf{T}_g$, and we define the {\em moduli space} ${\mathcal M}_g$ to be the quotient space ${\mathbf T}_g/{\rm Mod}_g$, see~\cite{Bro90, Har71, MS}. The projection $\mathbf{T}_g\rightarrow{\mathcal M}_g$ is a regular branched covering, so ${\mathcal M}_g$ has the structure of an orbifold. For $g\ge 3$ the (orbifold-)singular locus ({\em branch locus of the covering)} $\mathcal{B}_g$ of ${\mathcal M}_g$ is formed by the Riemann surfaces with non-trivial automorphisms, whereas for $g=2$ it consists of those with other automorphisms in addition to the identity and the hyperelliptic involution. See~\cite{MS}, for example.

More generally, if $\Gamma$ is an abstract group with signature $\sigma$, then the Teichm\"uller space $\mathbb{T}(\Gamma)$ is the space of embeddings $r: \Gamma \to PSL_2(\mathbb{R})$, with $r(\Gamma)$ discrete in $PSL_2(\mathbb{R})$, modulo conjugation in $PSL_2(\mathbb{R})$; if $\Gamma$ has signature $\sigma = (\gamma; m_1,\dots, m_k)$ then $\mathbb{T}(\Gamma)$ is homeomorphic to a ball of dimension $d=6\gamma - 6 + 2k$ (see~\cite{Har71, MS, Sin72}). For example, if $\gamma = 0$ and $k = 3$ then $d = 0$ and $\mathbb{T}(\Gamma)$ is a point, giving a single conjugacy class of triangle groups of type $(m_1,m_2,m_3)$ in $PSL_2(\mathbb{R})$. Since ${\mathbb T}(\Gamma)$ depends only on $\sigma$ we can write it as $\mathbb{T}(\sigma)$. 
The modular group of $\Gamma$ is the quotient ${\rm Mod}(\Gamma):={\rm Aut}^+(\Gamma)/{\rm Inn}(\Gamma)$, and the moduli space of $\Gamma$ is the quotient ${\mathcal M}(\Gamma):={\mathbb T}(\Gamma)/{\rm Mod}(\Gamma)$. Any inclusion $\alpha:\Gamma\rightarrow\Gamma^{\prime}$ of Fuchsian groups induces an embedding ${\mathbb T}(\alpha):{\mathbb T}(\Gamma')\rightarrow {\mathbb T}(\Gamma)$ defined by $[r]\mapsto\lbrack r\circ \alpha]$ (see~\cite{Har71, MS, Sin72}). Any action $\psi$ of a finite group $G$ on ${\mathcal S}={\mathbb H}/K$ is determined by an inclusion $\alpha:K\to\Gamma$ (via the monodromy $\theta:\Gamma\to G$). Then (\cite{Bro90, Har71, MS})
\begin{equation}\label{stratif}
{\mathcal B}_g = \bigcup_{G, \psi} \overline{\mathcal M}_g^{\,G, \psi}
\end{equation}
where $\overline{\mathcal M}_g^{\,G, \psi}$ is a closed, irreducible algebraic subvariety of ${\mathcal M}_g$ defined as consisting of those Riemann surfaces $\mathcal S$ with a group of automorphisms conjugate to $G$ in ${\rm Mod}_g$ (the conjugacy class determined by $\psi$). Its interior ${\mathcal M}_g^{\,G, \psi}$, if non-empty, is a smooth, locally closed algebraic subvariety of ${\mathcal M}_g$, dense in $\overline{\mathcal M}_g^{\,G, \psi}$; by definition, it consists of those surfaces with full automorphism group conjugate to $G$ in ${\rm Mod}_g$, and is called an {\em equisymmetric stratum}.

Observe that ${\mathcal M}_g^{\,G, \psi}$ is empty if and only if the action $\psi$ of $G$ extends for each Riemann surface admitting the action $\psi$ (see~\cite{CI} for example). An action of $G$, with surface epimorphism $\theta:\Gamma\to G$ and $\ker\theta=K$, is said to extend to an action of a group $G'\ge G$ if and only if there is an abstract Fuchsian group $\Gamma'\ge\Gamma$ with a surface epimorphism $\theta':\Gamma'\to G'$ such that $\theta'|_{\Gamma}=\theta$, $\ker\theta'=K$, and such that ${\mathbb T}(\Gamma)$ and ${\mathbb T}(\Gamma')$ have the same dimension. In~\cite{Sin72} Singerman determined all pairs of signatures $(\sigma(\Gamma), \sigma(\Gamma'))$ for which it is possible to have an extension in the sense introduced here. The action is called {\em maximal} if it has no such extension with $G'\ne G$. 


\subsection{Jacobian varieties with a group action.}\label{Jacob}

Let $\mathcal S$ be a compact Riemann surface of genus $g\ge 2$. We denote by ${\mathcal H}^1({\mathcal S};{\mathbb C})$ the $g$-dimensional complex vector space of $1$-forms on $\mathcal S$, and by $H_1({\mathcal S};{\mathbb Z})$ the first integral homology group of $\mathcal S$. Recall that, as mentioned in the Introduction, the Jacobian variety
\[J{\mathcal S}:={\mathcal H}^1({\mathcal S};{\mathbb C})^*/H_1({\mathcal S};{\mathbb Z})\]
of $\mathcal S$ is an irreducible principally polarized abelian variety of dimension $g$. See, for example, \cite{BL, FK}. 

Let $\mathcal S$ admit an action of a finite group $G$. A collection $\{H_1,\ldots, H_t\}$ of distinct subgroups of $G$ is termed $G$-{\sl admissible\/} if
\[d_V^{H_1}+\cdots+d_V^{H_t}\le d_V\]
for each non-trivial irreducible complex representation $V$ of $G$, where $d_V^{H_j}$ is the dimension of the subspace of $V$ fixed under the action of $H_j$ and $d_V$ is the degree of $V$. This collection is termed {\sl admissible\/} if it is $G$-admissible for some group $G$. If $\{H_1,\ldots, H_t\}$ is admissible then, by~\cite[Theorem~1.1]{RR19}, $J\mathcal S$ decomposes, up to isogeny, as
\[J{\mathcal S}\sim J({\mathcal S}/H_1)\times\cdots\times J({\mathcal S}/H_t)\times A\]
for some abelian subvariety $A$ of $J\mathcal S$.
 

\section{Signatures}\label{sign}

If $G$ acts on a compact Riemann surface $\mathcal S$ of genus $g\ge 2$, with $G\cong\Gamma/K$ and with a signature
\[\sigma=(\gamma;m_1,\ldots, m_k)\]
as before, then the Riemann--Hurwitz formula, applied to the inclusion $K\le\Gamma$, states that
\begin{equation}\label{RH}
2(g-1)=|G|\left(2\gamma-2+\sum_{i=1}^k\left(1-\frac{1}{m_i}\right)\right).
\end{equation}
The ratio $\rho=\rho_{\sigma}:=|G|/(g-1)$ depends only on $\sigma$; with this notation, equation~(\ref{RH}) becomes
\begin{equation}\label{2/rhoeqn}
\frac{2}{\rho}=2\gamma-2+\sum_{i=1}^k\left(1-\frac{1}{m_i}\right).
\end{equation}

The following result, presumably well-known, will be useful later on. Although we will not explicitly cite it, we rely on both the result and the method of proof when compiling lists of signatures satisfying various conditions.

\begin{lemma}\label{RHlemma}
Given any rational $\alpha>0$ there are only finitely many sets of integers $\gamma\ge 0$ and $m_1,\ldots,m_k\ge 2$ satisfying the equation
\begin{equation}\label{alphaeqn}
\alpha=2\gamma-2+\sum_{i=1}^k\left(1-\frac{1}{m_i}\right).
\end{equation}
\end{lemma}

\noindent{\sl Proof.} Since $1-\frac{1}{m_i}\ge \frac{1}{2}$ for each $i$ we have
\[2\gamma+\frac{k}{2}\le \alpha+2,\]
so there are only finitely many possibilities for $\gamma$ and $k$. For any given $\gamma$ and $k$ we need to solve an equation of the form
\begin{equation}\label{betaeqn}
\beta=\sum_{i=1}^k\frac{1}{m_i}
\end{equation}
with a fixed $\beta>0$. We will use induction on $k$ to show that equation~(\ref{betaeqn}) has only finitely many solutions. If $k=1$ this result is trivial, so suppose that $k\ge 2$ and we have proved it for sums of $k-1$ terms. We may number the terms $m_i$ so that $m_1\ge\ldots\ge m_k$ (temporarily departing from our usual convention), in which case $m_k\le k/\beta$ so that there are only finitely many possibilities for $m_k$. For each of these values of $m_k$ we are solving an equation
\[\beta'=\sum_{i=1}^{k-1}\frac{1}{m_i}\]
with a fixed $\beta'>0$, and by the induction hypothesis this has only finitely many solutions. \hfill$\square$

\medskip

Despite this result, there is no uniform bound on the number of solutions for equations of the form~(\ref{alphaeqn}). For example, given any solution for $\alpha$ with $\gamma\ge 1$, one can create another by replacing $\gamma$ with $\gamma-1$ and in compensation  increasing $k$ by adjoining four terms $m_i=2$. One can iterate this substitution as often as required if $\gamma$ is large enough. Nevertheless, in practical applications to Riemann surfaces the number of solutions is usually rather small, as we will see later.

\begin{cor}\label{finrho}
Each $\rho>0$ corresponds to only finitely many cocompact signatures $\sigma$.  \hfill$\square$
\end{cor}

As is well-known, the smallest positive value of the right-hand side of equation~(\ref{2/rhoeqn}) is $1/42$, attained only by $\sigma=(2,3,7)$ and leading to the Hurwitz bound $\rho\le 84$. This, together with Corollary~\ref{finrho}, implies that there are only finitely many signatures $\sigma$ such that $\rho_{\sigma}\in{\mathbb Z}$. For each integer $\rho$ all such signatures $\sigma$ can be found by simple (if tedious) arithmetic. Those for integers $\rho\ge 8$ can be found in the Appendix of~\cite{BJ}. Adjoining those for $\rho=4, 5, 6$ and $7$ leads to the following list, which gives all $\sigma$ corresponding to integers $\rho\ge 4$. (The rather long lists of signatures for $\rho=1, 2$ and $3$ are omitted here since for such $\rho$, any group of order $\rho p$ is cyclic or  isomorphic to $G_{p,\rho}$, and these can be dealt with more easily by a different method.)

\begin{itemize}
\item $\rho=84$, $\sigma=(2,3,7)$;
\item $\rho=48$, $\sigma=(2,3,8)$;
\item $\rho=40$, $\sigma=(2,4,5)$;
\item $\rho=36$, $\sigma=(2,3,9)$;
\item $\rho=30$, $\sigma=(2,3,10)$;
\item $\rho=24$, $\sigma=(2,3,12)$, $(2,4,6)$, $(3,3,4)$;
\item $\rho=21$, $\sigma=(2,3,14)$;
\item $\rho=20$, $\sigma=(2,3,15)$, $(2,5,5)$;
\item $\rho=18$, $\sigma=(2,3,18)$;
\item $\rho=16$, $\sigma=(2,3,24)$, $(2,4,8)$;
\item $\rho=15$, $\sigma=(2,3,30)$, $(2,5,6)$, $(3,3,5)$;
\item $\rho=14$, $\sigma=(2,3,42)$;
\item $\rho=13$, $\sigma=(2,3,78)$;
\item $\rho=12$, $\sigma=(2,4,12)$, $(2,6,6)$, $(3,3,6)$, $(3,4,4)$, $(2,2,2,3)$;
\item $\rho=10$, $\sigma=(2,4,20)$, $(2,5,10)$;
\item $\rho=9$, $\sigma=(2,4,36)$, $(2,6,9)$, $(3,3,9)$;
\item $\rho=8$, $\sigma=(2,5,20)$, $(2,6,12)$, $(2,8,8)$, $(3,3,12)$, $(3,4,6)$, $(2,2,2,4)$;
\item $\rho=7$, $\sigma=(2,5,70)$, $(2,6,21)$, $(2,7,14), (3,3,21)$;
\item $\rho=6$, $\sigma=(2,7,42)$, $(2,8,24)$, $(2,9,18), (2,10,15), (2,12,12), (3,4,12), (3,6,6), (4,4,6)$, $(2,2,2,6), (2,2,3,3)$;
\item $\rho=5$, $\sigma=(2,11,110), (2,12,60), (2,14,35), (2,15,30), (2,20,20)$, $(3,4,60), (3,5,15), (3,6,10)$, $(4,4,10)$, $(5,5,5)$, $(2,2,2,10)$;
\item $\rho=4$, $\sigma=(3,7,42), (3,8,24), (3,9,18), (3,10,15), (3,12,12)$,  $(4,5,20), (4,6,12), (4, 8, 8)$, $(5,5,10), (6,6,6)$, $(2,2,3,6), (2,2,4,4), (2,3,3,3)$, $(2^{[5]})$, $(1;2)$.
\end{itemize}

For future reference, let us define
\[\Sigma:=\{\sigma\mid\rho_{\sigma}\in{\mathbb Z},\, \rho\ge 4\},\]
the set of all signatures in the above list. 


\section{Normal structure}\label{normal}

From now on we will assume that $G$ acts on some compact Riemann surface $\mathcal S$ of genus $g=p+1$ with $p$ prime, and that $\rho:=|G|/(g-1)\in{\mathbb Z}$, or equivalently, $p$ divides $|G|$. The following lemma describes the normal structure shared by almost all of the groups $G$ we shall study.

\begin{lemma}\label{Sylowlemma}
{\rm (a)} If $\rho$ is coprime to $p$ and has no divisor $d\ne 1$ such that $d\equiv 1$ mod~$(p)$ (thus in particular if $p>\rho$) then $G$ is a semidirect product $P\rtimes Q$ where $P\cong C_p $ and $Q$ has order $\rho$.

\smallskip

\noindent{\rm(b)} If, in addition, $p$ is coprime to all the elliptic periods in the signature $\sigma$ of $G$, then $Q$ has a faithful action, with signature $\sigma$, as a group of automorphisms of the Riemann surface ${\mathcal T}:={\mathcal S}/P$ of genus $2$.
\end{lemma}

\noindent{\sl Proof.} (a) Since $|G|=\rho p$ is divisible by $p$ but not by $p^2$, $G$ has a Sylow $p$-subgroup $P\cong C_p$. Sylow's theorems state that the number of Sylow $p$-subgroups divides $|G|$ and is congruent to $1$ mod~$(p)$, so it divides $\rho$; by our hypothesis on the divisiors of $\rho$ this number must be $1$, so that $P$ is a normal subgroup of $G$. Since $|G:P|=\rho$ is coprime to $p$ the Schur--Zassenhaus Theorem implies that $G$ is a semidirect product $P\rtimes Q$ for some subgroup $Q$ of order $\rho$.

(b) The normal subgroup $P$ of $G$ lifts back, under the epimorphism $\theta:\Gamma\to G$ with kernel $K\cong\Pi_g$, to a normal subgroup $\Delta$ of $\Gamma$ with $\Gamma/\Delta\cong G/P\cong Q$ and $\Delta/K\cong C_p$. If $p$ does not divide any elliptic period of $\Gamma$ then $\Delta$ is torsion-free and therefore a surface group. Since the $p$-sheeted covering ${\mathcal S}={\mathbb H}/K\to{\mathbb H}/\Delta={\mathcal T}:={\mathcal S}/P$ is smooth, ${\mathcal T}$ has genus
\[\frac{g-1}{p}+1=2.\]
Since $Q\cong\Gamma/\Delta$ it follows that $Q$ acts faithfully as a group of automorphisms of ${\mathcal T}$. The obvious composition $\Gamma\to G\to Q$ is a surface epimorphism, so this action of $Q$ has the same signature as that of $G$, namely the signature $\sigma$ of $\Gamma$. \hfill$\square$

\medskip

(Note that if $g-1$ is a prime-power $p^e$, all of this lemma remains valid apart from the isomorphism of $P$ with $C_p$; this suggests an obvious generalisation of the present investigation, as in~\cite{CRC} for example.)

\begin{figure}[h!]
\begin{center}
\begin{tikzpicture}[scale=0.5, inner sep=0.8mm]

\node at (-10,-3) {$\Gamma$};
\node at (-10,-6) {$\Delta$};
\node at (-10,-9) {$K$};
\node at (-10,-12) {$\Delta'\Delta^p$};
\node at (-10,-15) {$1$};

\node at (-5,-6) {$M$};
\node at (-5,-9) {$\overline K$};
\node at (-5,-12) {$0$};

\node at (0,-3) {$G$};
\node at (0,-6) {$P$};
\node at (0,-9) {$1$};

\node at (5,-3) {$Q$};
\node at (5,-6) {$1$};

\node at (-7.5,-6) {$\to$};
\node at (-7.5,-9) {$\to$};
\node at (-7.5,-12) {$\to$};
\node at (-5,-3) {$\to$};
\node at (-2.5,-6) {$\to$};
\node at (-2.5,-9) {$\to$};

\node at (2.5,-3) {$\to$};
\node at (2.5,-6) {$\to$};

\draw [thick] (-10,-3.7) to (-10,-5.4);
\draw [thick] (-10,-6.7) to (-10,-8.4);
\draw [thick] (-10,-9.7) to (-10,-11.4);
\draw [thick] (-10,-12.7) to (-10,-14.4);
\draw [thick] (-5,-6.7) to (-5,-8.4);
\draw [thick] (-5,-9.7) to (-5,-11.4);
\draw [thick] (0,-3.7) to (0,-5.4);
\draw [thick] (0,-6.7) to (0,-8.4);
\draw [thick] (5,-3.7) to (5,-5.4);

\end{tikzpicture}

\end{center}
\caption{Normal structure of $\Gamma$ and $G$ in Lemma~\ref{Sylowlemma}}
\label{normal}
\end{figure}

The groups of automorphisms of Riemann surfaces of genus $2$, together with their corresponding signatures $\sigma$, are listed by Broughton in~\cite[Table 4]{Bro}. For each signature $\sigma$, all but finitely many primes $p$ satisfy the conditions of Lemma~\ref{Sylowlemma}. For such primes, since $P$ is abelian and of exponent $p$, $K$ must contain the commutator subgroup $\Delta'$ of $\Delta$ and the group $\Delta^p$ generated by its $p$th powers. Thus $\Delta\ge K\ge \Delta'\Delta^p$, so that $K$ projects onto a codimension $1$ submodule $\overline K:=K/\Delta'\Delta^p$ of the ${\mathbb F}_pQ$-module $M:=\Delta/\Delta'\Delta^p$. Figure~\ref{normal} shows the normal structure of $\Gamma$ and $G$ in Lemma~\ref{Sylowlemma}; the vertical lines denote inclusions of subgroups or submodules, and the arrows denote natural epimorphisms.

By its definition, $M$ is isomorphic, as an ${\mathbb F}_pQ$-module, to the mod~$(p)$ homology group $H_1({\mathcal T};{\mathbb F}_p)$ of $\mathcal T$. By decomposing this homology representation of $Q$ one can determine those primes $p$ which give examples of the actions we require. The small number of exceptional primes, for which Lemma~\ref{Sylowlemma} does not apply, can be dealt with individually.

Since $\mathcal T$ has genus $2$, $M$ has dimension $4$ over ${\mathbb F}_p$. Since $p$ does not divide $\rho=|Q|$, Maschke's Theorem applies to the action of $Q$ on $M$, so $M$ is a direct sum of irreducible submodules. Now $H_1({\mathcal T};{\mathbb C})=H_1({\mathcal T};{\mathbb Z})\otimes{\mathbb C}$ is a direct sum of two $Q$-submodules, corresponding under duality to the holomorphic and antiholomorphic differentials in $H^1({\mathcal T};{\mathbb C})$ and affording complex conjugate representations of $Q$~\cite{Sah}; this implies that $M$ is either irreducible, or a direct sum of two irreducible $2$-dimensional submodules or four irreducible $1$-dimensional submodules. Thus $M$ has a $1$-dimensional quotient if and only if the last case arises, giving four kernels $K$ corresponding to two chiral pairs of surfaces $\mathcal S$ of genus $g$. A theorem of Serre (see~\cite[V.3.4]{FK}, for example) shows that $Q$ acts faithfully on $H_1({\mathcal T};{\mathbb F}_p)$ for $p>2$, so this action embeds $Q$ in $GL_1(p)^4\cong C_{p-1}^4$, and hence $Q$ is an abelian group of rank at most $4$ and exponent $e$ dividing $p-1$.

The only abelian groups $Q$ in Broughton's list~\cite[Table~4]{Bro} of genus~$2$ group actions are the following:
\begin{enumerate}
\item $C_6\times C_2$ with $\sigma=(2,6,6)$;
\item $C_{10}$ with $\sigma=(2,5,10)$;
\item $C_8$ with $\sigma=(2,8,8)$;
\item $C_6$ with $\sigma=(3,6,6)$;
\item $C_6$ with $\sigma=(2,2,3,3)$;
\item $C_5$ with $\sigma=(5,5,5)$;
\item $C_4$ with $\sigma=(2,2,4,4)$;
\item $V_4$ with $\sigma=(2^{[5]})$;
\item $C_3$ with $\sigma=(3^{[4]})$;
\item $C_2$ with $\sigma=(2^{[6]})$;
\item $C_2$ with $\sigma=(1;2,2)$.

\end{enumerate}
Since we are restricting attention to integers $\rho\ge 3$, only cases~(1) to (9) are relevant here. (As we will see later, they correspond to cases~(i) to (ix) respectively in Theorem~\ref{integerthm}(a).) In each case one can use character theory to decompose the module $M$ for different primes $p$, and thus determine those giving $1$-dimensional quotients and how $Q$ acts on them.

By the Lefschetz fixed-point formula, the homology character $\chi$ of $Q$ on $H_1({\mathcal S};{\mathbb Z})$ is $2-\phi$ where $\phi(q)$ is the number of fixed points of an element $q\in Q$ on $\mathcal T$. By a result of Macbeath~\cite{Macb73}, 
\begin{equation}\label{Macbeatheqn}
\phi(q)=|N_Q(\langle q\rangle)|\sum_{i=1}^k\frac{\varepsilon_i(q)}{m_i}
\end{equation}
for all $q\ne 1$ in $Q$, where $\varepsilon_i(q)=1$ or $0$ as $q$ is or is not conjugate in $Q$ to a power of the image of the elliptic generator $X_i$ of $\Gamma$. When $Q$ is cyclic this simplifies to
\[\phi(q)=|Q|\sum\frac{1}{m_i}, \]
where the sum is over all $m_i$ divisible by the order of $q$. Using the resulting values of $\chi$, one can calculate the coefficients
\[a_j=\frac{1}{|Q|}\sum_{q\in Q}\chi(q)\chi_j(q)\]
of the irreducible characters $\chi_j$ of $Q$ in the decomposition
\[\chi=\sum_j a_j\chi_j\]
of $\chi$ as a sum of irreducible complex characters of $Q$. Since $p\equiv 1$ mod~$(e)$, ${\mathbb F}_p$ is a splitting field for $Q$, so reducing this decomposition mod~$(p)$ gives the decomposition of $M$ over ${\mathbb F}_p$. 

An alternative method of evaluating $\phi$ is to use an explicit model for $\mathcal T$ and $Q$, and simply to find and count the fixed points of each $q\in Q$, as in~\cite{BJ}. For instance, in case (1) one can take $\mathcal T$ to be the compact Riemann surface corresponding to the curve $w^2=z^6-1$, with the direct factors of $G$ generated by its automorphisms $z\mapsto e^{\pi i/3}z$ and $w\mapsto -w$; in cases~(2) and (3) one can use the curves $w^2=z^5-1$ and $w^2=z(z^4-1)$ with $G$ generated by $(z,w)\mapsto(e^{2\pi i/5}z,-w)$ and $(z,w)\mapsto (iz,e^{\pi i/4}w)$ respectively, and in cases~(4) to (9) one can restrict $\phi$ to subgroups of these three groups.


\section{Proof of Theorem~\ref{integerthm} for good primes}

Here we will prove parts (a) and (b) of Theorem~\ref{integerthm} for `good primes' $p$, those satisfying the conditions of Lemma~\ref{Sylowlemma}, by considering the possibilities for $Q$ in cases~(1) to (9) in turn. The remaining `bad primes' will be dealt with in the next section. We  will also prove Theorem~\ref{integerthm}(c) for some of cases~(i) to (xii) in this section, postponing others until parts (a) and (b) have been proved.

The Teichm\"uller space ${\mathbb T}(\sigma)$ of groups $\Gamma$ of a given signature $\sigma=(\gamma;m_1,\ldots, m_k)$ has dimension $d=6\gamma-6+2k$ (see~\cite{Sin72}, for example). In cases (1), (2), (3), (4) and (6), where $\gamma=0$ and $k=3$, we have $d=0$ and ${\mathbb T}(\sigma)$ is a point, giving a single conjugacy class of triangle groups $\Gamma$ in $PSL_2({\mathbb R})$. However, in the remaining cases $d>0$ and we have $d$-dimensional families of groups $\Gamma$ and of surfaces $\mathcal S$. We will deal with the triangle groups first, since the decomposition of the homology character in these cases has already been determined by Kazaz~\cite {Kaz} in the context of hypermaps of genus~$2$ and their coverings; the other cases follow easily.

We will deal with case~(1) in some detail, and then just outline the method and results for the other cases. Here $\Gamma=\Gamma(2,6,6)$, $\Delta=\Gamma'$ and $Q\cong C_6\times C_2\cong V_4\times C_3$ of order $\rho=12$ and exponent $e=6$.  In this case all primes $p>5$ satisfy the conditions of Lemma~\ref{Sylowlemma} (note that for $p=5$ a group $G$ of order $\rho p=60$, such as $A_5$, could have six Sylow $5$-subgroups $P$, rather than one). Let $x_1$ and $x_2$ be the images in $Q$ of the elliptic generators $X_1$ and $X_2$ of $\Gamma$ (see \S\ref{sign}), generating direct factors $C_2$ and $C_6$ of $Q$. The irreducible complex characters of $Q$ are the homomorphisms $\chi_{i,j}=\chi_1^i\chi_2^j:Q\to S^1$ ($i\in{\mathbb Z}_2, j\in{\mathbb Z}_6$), where $\chi_1:x_1\mapsto -1, x_2\mapsto 1$ and $\chi_2:x_1\mapsto 1, x_2\mapsto\zeta$ (a primitive 6th root of $1$). Using equation~(\ref{Macbeatheqn}) we find that $\chi=\chi_{1,1}+\chi_{1,-1}+\chi_{1,2}+\chi_{1,-2}$. The first two and the last two are complex conjugate pairs, with image $C_6$ and kernels generated by the involutions $x_2^3$ and $x_1x_2^3$ respectively. It follows that for primes $p\equiv 1$ mod~$(6)$ (equivalently, $p\equiv 1$ mod~$(3)$), $Q$ acts as $C_6\times C_2$ on each of four corresponding submodules $M_{i,j}$ of $M$, with the factors $C_6$ acting faithfully and $C_2$ trivially. It acts in the same way on the corresponding $1$-dimensional quotient modules of $M$, each obtained by factoring out the other three submodules, so it induces groups $G\cong G_{p,6}\times C_2$ on two chiral pairs of surfaces $\mathcal S$ of genus $p+1$.

Now the normaliser $N(\Gamma)$ of $\Gamma$ in $PSL_2({\mathbb R})$ is the maximal Fuchsian group $\Gamma(2,4,6)$, which contains $\Gamma$ with index $2$ (see~\cite{Sin72}). Conjugation in $N(\Gamma)$ transposes the elliptic generators $X_2$ and $X_3=(X_1X_2)^{-1}$ of $\Gamma$, so it acts on $Q$ by transposing the involutions $x_2^3$ and $x_1x_2^3$. It therefore transposes the first chiral pair of surfaces with the second, so up to isomorphism we have just one chiral pair of surfaces $\mathcal S_1$ and $\overline{\mathcal S}_1$. It also follows from~\cite{Sin72} that $N(\Gamma)$ is the only Fuchsian group properly containing $\Gamma$, so each surface has automorphism group $G$. This deals with case~(1), giving Theorem~\ref{integerthm}(a)(i) together with the statements in parts (b) and (c) concerning this case.

\begin{figure}[h!]
\begin{center}
\begin{tikzpicture}[scale=0.5, inner sep=0.8mm]

\node at (-23,3) {(i)};
\node at (-23,0) {(v), (iv)};
\node at (-23,-3) {(ix)};

\node at (-12.5,6) {$N(\Gamma)=(2,4,6)$};
\node at (-12.5,3) {$\Gamma=(2,6,6)$};
\node at (-18.2,0) {$\Gamma_1=(2,2,3,3)$};
\node at (-12.5,0) {$\Gamma_2=(3,6,6)$};
\node at (-9.7,0) {$\sim$};
\node at (-6.8,0) {$\Gamma_3=(3,6,6)$};
\node at (-12.5,-3) {$\Gamma_4=(3^{[4]})$};
\node at (-12.5,-6) {$\Delta$};

\node at (0,3) {$G=G_{p,6}\times C_2$};
\node at (-2,0) {$G_1$};
\node at (0,0) {$G_2$};
\node at (2,0) {$G_3$};
\node at (0,-3) {$G_{p,3}$};
\node at (0,-6) {$P$};

\node at (6.5,3) {$Q\cong C_6\times C_2$};
\node at (6.5,0) {$Q_i\cong C_6$};
\node at (6.5,-3) {$C_3$};
\node at (6.5,-6) {$1$};

\node at (-3.5,3) {$\to$};
\node at (-3.5,0) {$\to$};
\node at (-3.5,-3) {$\to$};
\node at (-3.5,-6) {$\to$};

\node at (3.5,3) {$\to$};
\node at (3.5,0) {$\to$};
\node at (3.5,-3) {$\to$};
\node at (3.5,-6) {$\to$};

\draw [thick] (-12.5,5.5) to (-12.5,3.8);
\draw [thick] (-13.5,2.5) to (-16.5,0.8);
\draw [thick] (-12.5,2.5) to (-12.5,0.8);
\draw [thick] (-11.5,2.5) to (-8.5,0.8);
\draw [thick] (-16.5,-0.7) to (-13.5,-2.4);
\draw [thick] (-12.5,-0.7) to (-12.5,-2.4);
\draw [thick] (-8.5,-0.7) to (-11.5,-2.4);
\draw [thick] (-12.5,-3.7) to (-12.5,-5.4);

\draw [thick] (-0.5,2.5) to (-1.5,0.8);
\draw [thick] (0,2.5) to (0,0.8);
\draw [thick] (0.5,2.5) to (1.5,0.8);
\draw [thick] (-1.5,-0.7) to (-0.5,-2.4);
\draw [thick] (0,-0.7) to (0,-2.4);
\draw [thick] (1.5,-0.7) to (0.5,-2.4);
\draw [thick] (0,-3.7) to (0,-5.4);

\draw [thick] (-11.5,2.5) to (-8.5,0.8);
\draw [thick] (-16.5,-0.7) to (-13.5,-2.4);
\draw [thick] (-12.5,-0.7) to (-12.5,-2.4);
\draw [thick] (-8.5,-0.7) to (-11.5,-2.4);

\draw [thick] (6.5,2.5) to (6.5,0.8);
\draw [thick] (6.5,-0.7) to (6.5,-2.4);
\draw [thick] (6.5,-3.7) to (6.5,-5.4);

\end{tikzpicture}

\end{center}
\caption{Cases (i), (iv), (v) and (ix) of Theorem~\ref{integerthm}(a)}
\label{1459}
\end{figure}
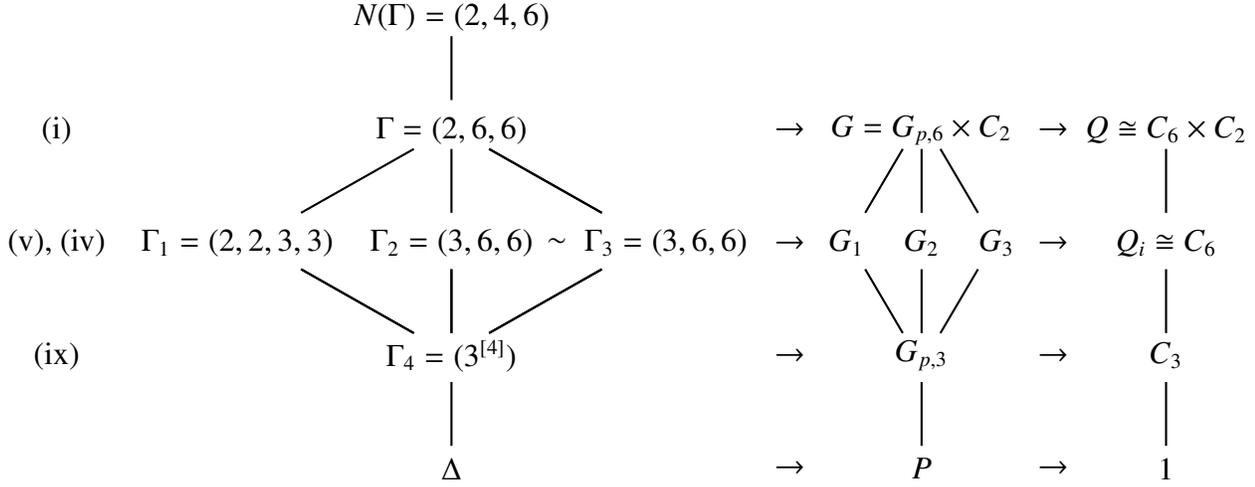

We can also deal with case~(4), and gain some information about cases~(5), (8) and (9), by considering subgroups of $\Gamma(2,6,6)$. This group has three subgroups of index $2$: the normal closure $\Gamma_1$ of $X_1$ and $X_2^2$ has signature $(2,2,3,3)$, while the normal closures $\Gamma_2$ and $\Gamma_3$ of $X_2$ and $X_3$ have signature $(3,6,6)$. These correspond to three subgroups $Q_i\cong C_6$ of index $2$ in the group $Q\cong C_6\times C_2$ in case~(1), and to three subgroups $G_i$ of index $2$ in $G\cong C_{p,6}\times C_2$, each acting on the same chiral pair of surfaces $\mathcal S_1$ and $\overline{\mathcal S}_1$ as in Theorem~\ref{integerthm}(a)(i). These subgroups are shown in Figure~\ref{1459}, where to save space we have represented a Fuchsian group $\Gamma$ by its signature $\sigma$, and the three subgroups $Q_1, Q_2$ and $Q_3$ by a single symbol $Q_i$.

In dealing with case~(4), since there is only one conjugacy class of subgroups $\Gamma(3,6,6)$ in $PSL_2({\mathbb R})$, we may without loss of generality take $\Gamma=\Gamma_2$ or $\Gamma_3$; indeed, since these are conjugate in $\Gamma(2,4,6)$ it is sufficient to consider just one of them, say $\Gamma_2$. As in case~(1), $\Gamma_2$ has a unique normal surface subgroup with the required quotient (now $C_6$), so this must be the same subgroup $\Delta$ as in case~(1), with the same quotient $M=\Delta/\Delta'\Delta^p$, but now regarded as a module for ${\mathbb F}_pQ_2$. This decomposes in the same way as before, so we obtain the same two chiral pairs of surfaces, isomorphic under $\Gamma(2,4,6)$, and we need to determine how $\Gamma_2$ acts on them. Now $Q_2=\langle x_2\rangle$, which contains the kernel $\langle x_2^3\rangle$ of $\chi_{1,\pm 1}$ but not the kernel $\langle x_1x_2^3\rangle$ of $\chi_{1,\pm 2}$, so $\Gamma_2$ induces a group $G_{p,3}\times C_2$ on the first chiral pair of surfaces and $G_{p,6}$ on the second; for $\Gamma_3$ it is the other way round. Thus $\mathcal S_1$ and $\overline{\mathcal S}_1$ each admit two groups $G_{p,3}\times C_2$ and $G_{p,6}$, as described in Theorem~\ref{integerthm}(a)(iv).

The subgroup $Q_1=\langle x_1, x_2^2\rangle=\langle x_1x_2^2\rangle$ contains neither of the kernels $\langle x_2^3\rangle$ and $\langle x_1x_2^3\rangle$ of $\chi_{1,\pm 1}$ and $\chi_{1,\pm 2}$, so $\Gamma_1$ induces a group $G_1\cong G_{p,6}$ on each of $\mathcal S_1$ and $\overline{\mathcal S}_1$, corresponding to case~(5). However, $\Gamma_1$ is only one of a $2$-dimensional family of groups $\Gamma(2,2,3,3)$, each having a unique normal surface subgroup $\Delta$ with quotient $Q\cong C_6$. In this case, each of the two faithful $1$-dimensional characters of $Q$ has multiplicity $2$ in $\chi$, so $M$ is a direct sum of two copies each of two $1$-dimensional submodules affording the two faithful actions of $Q$. Now the direct sum of two isomorphic $1$-dimensional modules also contains another $p-1$ copies of that module, so $M$ has $2(p+1)$ maximal submodules, each corresponding to a kernel $K<\Delta$ and a surface $S$ affording a group $G=G_1\cong G_{p,6}$ (see Theorem~\ref{integerthm}(a)(v)).

The three subgroups $\Gamma_i \;(i=1,2,3)$ of index $2$ in $\Gamma=\Gamma(2,6,6)$ intersect in a normal subgroup $\Gamma_4$ of index $4$ and signature $(3^{[4]})$. This is one of a $2$-dimensional family of groups $\Gamma(3^{[4]})$ corresponding to the group $Q=C_3$ and signature $\sigma=(3^{[4]})$ in case~(9), giving groups $G\cong G_{p,3}$ for $p\equiv 1$ mod~$(3)$ as in Theorem~\ref{integerthm}(a)(ix) (the only other group of order $\rho(g-1)=3p$ for $p>3$, namely $C_{3p}$, is not generated by elements of order $3$). It is shown in~\cite[Corollary~1]{IR} that any action $\Gamma(3^{[4]})\to G_{p,3}$ extends to an action $\Gamma(2,2,3,3)\to G_{p,6}$ on the same family of surfaces as in Theorem~\ref{integerthm}(a)(v), where $\Gamma(2,2,3,3)$ and $G_{p,6}$ contain $\Gamma(3^{[4]})$ and $G_{p,3}$ with index~$2$ (note that the decomposition of the homology module is the same for $Q=C_3$ in case~(9) as for $Q=C_6$ in case~(5)). In fact, any group $\Gamma(3^{[4]})$ has three normal surface subgroups $\Delta$ of index $3$, equivalent under outer automorphisms and each yielding $2(p+1)$ maximal submodules of its quotient module $M=\Delta/\Delta'\Delta^p$, so we have $6(p+1)$ kernels $K$ in $\Gamma(3^{[4]})$.

\begin{figure}[h!]
\begin{center}
\begin{tikzpicture}[scale=0.5, inner sep=0.8mm]

\node at (-10,3) {(i)};
\node at (-10,0) {(viii)};

\node at (-6.5,3) {$\Gamma=(2,6,6)$};
\node at (-6.5,0) {$\Gamma_0=(2^{[5]})$};
\node at (-6.5,-3) {$\Delta$};

\node at (0,3) {$G=G_{p,6}\times C_2$};
\node at (0,0) {$G_{p,2}\times C_2$};
\node at (0,-3) {$P$};

\node at (6.5,3) {$Q\cong C_6\times C_2$};
\node at (6.5,0) {$Q_0\cong V_4$};
\node at (6.5,-3) {$1$};

\node at (-3.5,3) {$\to$};
\node at (-3.5,0) {$\to$};
\node at (-3.5,-3) {$\to$};

\node at (3.5,3) {$\to$};
\node at (3.5,0) {$\to$};
\node at (3.5,-3) {$\to$};

\draw [thick] (-6.5,2.5) to (-6.5,0.8);
\draw [thick] (-6.5,-0.7) to (-6.5,-2.4);

\draw [thick] (0,2.5) to (0,0.8);
\draw [thick] (0,-0.7) to (0,-2.4);

\draw [thick] (6.5,2.5) to (6.5,0.8);
\draw [thick] (6.5,-0.7) to (6.5,-2.4);

\end{tikzpicture}

\end{center}
\caption{Cases (i) and (viii) of Theorem~\ref{integerthm}(a)}
\label{1and8}
\end{figure}

There is also a unique normal subgroup $\Gamma_0$ of index $3$ in $\Gamma(2,6,6)$, namely the normal closure of $X_1$ and $X_2^3$, with signature $(2^{[5]})$ and with $Q_0=\Gamma_0/\Delta\cong V_4$, corresponding to case~(8); see Figure~\ref{1and8}. Again, for $p\equiv 1$ mod~$(3)$ the chiral pair $\mathcal S_1$ and $\overline{\mathcal S}_1$ in case~(1) admit a corresponding subgroup of $G_{p,6}\times C_2$, this time isomorphic to $G_{p,2}\times C_2\cong D_p\times C_2\cong D_{2p}$ as described in Theorem~\ref{integerthm}(a)(viii); once again, they are members of a family of surfaces, this time of dimension $4$, admitting this group but now without the restriction that $p\equiv 1$ mod~$(3)$.

\begin{figure}[h!]
\begin{center}
\begin{tikzpicture}[scale=0.5, inner sep=0.8mm]

\node at (-10,3) {(ii)};
\node at (-10,0) {(vi)};

\node at (-6,6) {$N(5,5,5)=(2,3,10)$};
\node at (-6,3) {$\Gamma=(2,5,10)$};
\node at (-6,0) {$(5,5,5)$};
\node at (-6,-3) {$\Delta$};

\node at (0,3) {$G=G_{p,10}$};
\node at (0,0) {$G_{p,5}$};
\node at (0,-3) {$P$};

\node at (5,3) {$Q\cong C_{10}$};
\node at (5,0) {$C_5$};
\node at (5,-3) {$1$};

\node at (-3,3) {$\to$};
\node at (-3,0) {$\to$};
\node at (-3,-3) {$\to$};

\node at (2.5,3) {$\to$};
\node at (2.5,0) {$\to$};
\node at (2.5,-3) {$\to$};

\draw [thick] (-6,5.5) to (-6,3.8);
\draw [thick] (-6,2.5) to (-6,0.8);
\draw [thick] (-6,-0.7) to (-6,-2.4);

\draw [thick] (0,2.5) to (0,0.8);
\draw [thick] (0,-0.7) to (0,-2.4);

\draw [thick] (5,2.5) to (5,0.8);
\draw [thick] (5,-0.7) to (5,-2.4);

\end{tikzpicture}

\end{center}
\caption{Cases (ii) and (vi) of Theorem~\ref{integerthm}(a)}
\label{2and6}
\end{figure}

A simpler case than case~(1) is case~(2), where $\Gamma=\Gamma(2,5,10)$,  $\Delta=\Gamma'$ and $\rho=e=10$; see Figure~\ref{2and6}. Here all primes $p>5$ satisfy the conditions of Lemma~\ref{Sylowlemma}, and character theory shows that $M$ splits as a sum of $1$-dimensional irreducible submodules if and only if $p\equiv 1$ mod~$(10)$, or equivalently $p\equiv 1$ mod~$(5)$. The resulting quotient modules $\overline K$ realise the four faithful $1$-dimensional representations of $Q$ over ${\mathbb F}_p$, so they correspond to two chiral pairs of surfaces $\mathcal S_2, \overline{\mathcal S}_2$ and $\mathcal S_2', \overline{\mathcal S'_2}$ acted on by groups $G\cong G_{p,10}$, as stated in Theorem~\ref{integerthm}(a)(ii). Since $N(\Gamma)=\Gamma$ (see~\cite{Sin72}), the four kernels $K\le\Delta$ are mutually non-conjugate in $PSL_2({\mathbb R})$, so the four surfaces $\mathcal S$ they uniformise are mutually non-isomorphic, with ${\rm Aut}\,{\mathcal S}=G$. The index $2$ inclusion $\Gamma(5,5,5)<\Gamma(2,5,10)$ shows that the group $\Delta$ in case~(6), where $\rho=e=5$, is the same as in case~(2); again, we obtain surfaces $\mathcal S$ if and only if $p\equiv 1$ mod~$(5)$, as claimed in Theorem~\ref{integerthm}(a)(vi). These are the chiral pairs $\mathcal S_2, \overline{\mathcal S}_2$ and $\mathcal S_2', \overline{\mathcal S'_2}$ in case~(2), each surface admitting an action of $G_{p,5}$ as a subgroup of its automorphism group $G_{p,10}$. In fact, $\Gamma(5,5,5)$ has three normal surface subgroups $\Delta$ of index $5$, each yielding four kernels $K$; however $N(\Gamma(5,5,5))$ is $\Gamma(2,3,10)$, which contains $\Gamma(2,5,10)$ with index $3$, and the quotient $\Gamma(2,3,10)/\Gamma(5,5,5)\cong S_3$ permutes these three subgroups $\Delta$ transitively, so up to isomorphism we obtain only the four surfaces described here. (The distinction between these twelve kernels becomes important when we consider them in Section~\ref{maps} as representing distinct hypermaps.)

\begin{figure}[h!]
\begin{center}
\begin{tikzpicture}[scale=0.5, inner sep=0.8mm]

\node at (-10,3) {(iii)};
\node at (-10,0) {(vii)};

\node at (-6,6) {$N(\Gamma)=(2,4,8)$};
\node at (-6,3) {$\Gamma=(2,8,8)$};
\node at (-6,0) {$(2,2,4,4)$};
\node at (-6,-3) {$\Delta$};

\node at (0,3) {$G=G_{p,8}$};
\node at (0,0) {$G_{p,4}$};
\node at (0,-3) {$P$};

\node at (5,3) {$Q\cong C_8$};
\node at (5,0) {$C_4$};
\node at (5,-3) {$1$};

\node at (-3,3) {$\to$};
\node at (-3,0) {$\to$};
\node at (-3,-3) {$\to$};

\node at (2.5,3) {$\to$};
\node at (2.5,0) {$\to$};
\node at (2.5,-3) {$\to$};

\draw [thick] (-6,5.5) to (-6,3.8);
\draw [thick] (-6,2.5) to (-6,0.8);
\draw [thick] (-6,-0.7) to (-6,-2.4);

\draw [thick] (0,2.5) to (0,0.8);
\draw [thick] (0,-0.7) to (0,-2.4);

\draw [thick] (5,2.5) to (5,0.8);
\draw [thick] (5,-0.7) to (5,-2.4);

\end{tikzpicture}

\end{center}
\caption{Cases (iii) and (vii) of Theorem~\ref{integerthm}(a)}
\label{3and7}
\end{figure}

The situation in case~(3), where $\rho=e=8$, is similar to that in case~(2); see Figure~\ref{3and7}. Now there are two normal subgroups of $\Gamma=\Gamma(2,8,8)$ with quotient $C_8$, but only one of them, the normal closure $\Delta$ of $X_1X_2^4$ in $\Gamma$, is torsion-free and thus a surface group of genus $2$. In this case Lemma~\ref{Sylowlemma} applies to all primes $p>7$, together with $p=5$, and $M$ has $1$-dimensional quotients if and only if $p\equiv 1$ mod~$(8)$. Again, these realise the four faithful $1$-dimensional representations of $Q$, so we obtain two chiral pairs of surfaces $\mathcal S$ admitting actions of $G\cong G_{p,8}$. This time, however, $\Gamma$ is not a maximal Fuchsian group: it has index $2$ in its normaliser $N(\Gamma)=\Gamma(2,4,8)$, which is maximal. Conjugation in $N(\Gamma)$, which leaves $\Delta$ invariant, induces isomorphisms between the two chiral pairs, so up to isomorphism we obtain one chiral pair $\mathcal S_3$ and $\overline{\mathcal S}_3$, as claimed in Theorem~\ref{integerthm}(a)(iii), each with ${\rm Aut}\,{\mathcal S}=G$. 

There is a unique subgroup of index $2$ in $\Gamma=\Gamma(2,8,8)$ containing $\Delta$; this has signature $(2,2,4,4)$ and $\rho=4$, corresponding to case~(7) where $Q\cong C_4$. This induces actions of $G_{p,4}$ on the same chiral pair of surfaces $\mathcal S_3$ and $\overline{\mathcal S}_3$ as in case~(3), where $p\equiv 1$ mod~$(8)$. However, as in case~(5), these are members of a $2$-dimensional family of groups $\Gamma(2,2,4,4)$ and surfaces $\mathcal S$ realising $G_{p,4}$, which arise for all primes $p\equiv 1$ mod~$(4)$ as in Theorem~\ref{integerthm}(a)(vii).

We have now dealt with cases~(1) to (9), corresponding to cases~(i) to (ix) in Theorem~\ref{integerthm}(a). These are the cases where the prime $p$ satisfies the conditions of Lemma~\ref{Sylowlemma}, so that $G$ has the normal structure $C_p\rtimes Q$ described there. Cases~(x) to (xii), where the lemma does not apply, will be considered in the next section.


\section{Exceptional actions for bad primes}

For each integer $\rho\ge 3$ and its corresponding signatures $\sigma\in\Sigma$ we need to consider the `bad primes' $p\ge 7$ which do not satisfy the conditions of Lemma~\ref{Sylowlemma}, since for these the arguments of the preceding section do not apply, and exceptional groups $G$ may appear. The primes dividing each $\rho$, or dividing $d-1$ for divisors $d\ne 1$ of $\rho$ (see Lemma~\ref{Sylowlemma}(a)), are easily found. For $\rho=84$ they are $2, 3, 5, 7, 11, 13, 41$ and $83$. For $\rho=48$ they are $2, 3, 5, 7, 11, 23$ and $47$. The primes $29$ and $23$ arise for $\rho=30$ and $24$ respectively, while $19$ arises for $\rho=40$ and $20$, and $17$ arises for $\rho=36$ and $18$. For other $\rho$, only primes $p\le 13$ arise. The largest primes dividing any of the elliptic periods in the signatures $\sigma\in\Sigma$ (see Lemma~\ref{Sylowlemma}(b)) are $13$ when $\rho=13$ and $\sigma=(2,3,78)$, and $11$ when $\rho=5$ and $\sigma=(2,11,110)$, so Lemma~\ref{Sylowlemma}(b) applies for all $p\ge 17$, and in some cases for smaller primes. 

In fact, a case-by-case argument in~\cite{BJ} shows that if $\rho\ge 8$ there are no exceptional groups for primes $p\ge 17$, although there is one for $p=13$, namely $PSL_2(13)$, which acts on three Riemann surfaces of genus $14$ as a Hurwitz group (with $\rho=84$) as in Theorem~\ref{integerthm}(a)(x). By inspection, the only bad primes arising for $3\le\rho\le 7$ are $2, 3, 5$ and $7$, so it is sufficient to restrict attention to the primes $p\le 13$. 

In addition to $\rho=84$ with $\sigma=(2,3,7)$, the prime $p=13$ is bad for $\rho=40$ with $\sigma=(2,4,5)$, for $\rho=14$ with $\sigma=(2,3,42)$, and $\rho=13$ with $\sigma=(2,3,78)$. There is no surface epimorphism from $\Gamma(\sigma)$ to a group $G$ of order $13\rho$ in the second or third of these three cases, since $|G|$ is coprime to $3$. In the first case, applying Sylow's theorems for the primes $5$ and $13$ shows that a group of order $13\rho=2^3.5.13$ has a normal subgroup of order $65$, whereas there is no epimorphism from $\Gamma(2,4,5)$ to the resulting quotient group of order $8$. Thus the only exceptional group $G$ arising for $p=13$ is the unique Hurwitz group $PSL_2(13)$ of genus $14$.

The prime $p=11$ is bad for all $\sigma$ with $\rho$ divisible by $12$, together with $\rho=5$ for $\sigma=(2,11, 110)$. The last case can be eliminated, since a group of order $55$ can have no element of order $110$. In the other cases $G$ acts by conjugation as a doubly transitive group of degree $12$ on its Sylow $11$-subgroups. Now the doubly transitive finite groups are known (see~\cite[Section 7.7]{DM}, for example), and those of degree $12$ are $PSL_2(11)$ and $PGL_2(11)$ acting naturally on ${\mathbb P}^1({\mathbb F}_{11})$, the Mathieu groups $M_{11}$ and $M_{12}$ acting on the cosets of a subgroup $PSL_2(11)$ and on the Steiner system $S(12, 6, 5)$, and $A_{12}$ and $S_{12}$ acting naturally; these groups all have orders divisible by $5$ whereas $G$ does not, so there are no exceptional groups for $p=11$.

The prime $p=7$ is bad for many signatures, including $(2,3,8)$ and $(3,3,4)$ for $\rho=48$ and $24$. These lead to two exceptional groups, namely $PGL_2(7)$ and its subgroup $PSL_2(7)$, as in Theorem~\ref{integerthm}(a)(xi) and (xii); each of these is (by character theory and M\"obius inversion, see~\cite[Sections 5.1.5, 5.1.6]{JW}, for example) a quotient of the corresponding triangle group by two normal subgroups $K$, so they both act on the same pair of surfaces of genus $8$. All other signatures can be eliminated by group-theoretic arguments as above, or (less laboriously) by checking Conder's lists of group actions~\cite{ConLarge} for examples of genus $8$ with $\rho\in{\mathbb Z}$. This completes the proof of Theorem~\ref{integerthm}(a) and (b).


\section{Proof of Theorem~\ref{integerthm}(c)}\label{partc}

In proving parts (a) and (b) of Theorem~\ref{integerthm}, we have described all the pairs $\mathcal S$ and $G\le A:={\rm Aut}\,{\mathcal S}$ for primes $p=g-1\ge 7$ and ratios $\rho=|G|/(g-1)\ge 3$. If $A\ne G$ then $A$ must be one of the groups described as acting on $\mathcal S$, but with ratio $\rho_A:=|A|/(g-1)$ a proper multiple of $\rho$. It is straightforward to check parts (a) and (b) to determine when this is possible, starting with the largest values of $\rho$.

Clearly, in cases (x) and (xi) no proper multiples of $\rho$ arise, so $A=G$ in these cases. (Indeed, the corresponding two triangle groups are maximal, giving a more direct proof.) In case (xii), however, $G=PSL_2(7)$ is a subgroup of index $2$ in the group $PGL_2(7)$ in case (xi), acting on the same two surfaces, so here $G<A=PGL_2(7)$.

The only proper multiples of $\rho=12$ in case~(i) are the values $84$, $48$ and $24$ in cases (x) to (xii). However, the groups $PSL_2(13)$, $PGL_2(7)$ and $PSL_2(7)$ in those three cases do not contain subgroups isomorphic to $G_{p,6}\times C_2$ for any prime $p$ (see~\cite{CCNPW} for their maximal subgroups, for example), so $A=G$ for all $p$ in case~(i). The surfaces in~(iv) are the same chiral pair ${\mathcal S}_1$ and $\overline{\mathcal S}_1$ as in~(i), so in case~(iv) we have $G<A\cong G_{p,6}\times C_2$ for each $p$.

No proper multiples of $\rho=10$ in case~(ii) appear, so $A=G$ in this case. This is the only proper multiple of $\rho=5$ in case~(vi), and we have seen that the surfaces in case~(ii) and (vi) are the same, so for case~(vi) we have $G<A\cong G_{p,10}$.

The only proper multiples of $\rho=8$ in case~(iii) are in cases~(xi) and (xii); however, the latter require $p=7$ whereas $p\equiv 1$ mod~$(8)$ in case~(iii), so here $A=G$.

The only proper multiples of $\rho=6$ in case~(v) are in cases~(i), (x), (xi) and (xii). A group $G_{p,6}$ in case~(v) cannot be a subgroup of $PSL_2(7)$, since this group has Sylow $7$-normalisers isomorphic to $G_{7,3}$, but it is a subgroup of $G_{p,6}\times C_2$, acting on ${\mathcal S}_1$ and $\overline{\mathcal S}_1$, for all $p\equiv 1$ mod~$(3)$, and of $PGL_2(7)$ and $PSL_2(13)$ for $p=7$ and $13$ respectively. Thus if ${\mathcal S}={\mathcal S}_1$ or $\overline{\mathcal S}_1$ then $G<A\cong G_{p,6}\times C_2$, and if $p=7$ or $13$ there are two or three surfaces $\mathcal S$ with $G<A\cong PSL_2(7)$ or $PSL_2(13)$, but otherwise $A=G$.

In case~(vii) there are proper multiples of $\rho=4$ in cases~(i), (iii), (x), (xi) and (xii). However, the requirements that $G_{p,4}\le A$ and $p\equiv 1$ mod~$(4)$ exclude all except (iii) (for example, the Sylow $13$-normaliser in $PSL_2(13)$ is isomorphic to $G_{13,6}$). We have seen that if $p\equiv 1$ mod~$(8)$ then $G_{p,4}$ acts as a subgroup of index~$2$ in the group $G_{p,8}$ in case~(ii) on the chiral pair ${\mathcal S}_3$ and $\overline{\mathcal S}_3$, so for these surfaces $G<A\cong G_{p,8}$, whereas $A=G$ for all other surfaces in case~(vii).

In case~(viii) there are also proper multiples of $\rho=4$ in cases~(i), (iii), (x), (xi) and (xii). The existence of elements of order $2p$ in $G$ excludes all except~(i), whereas if $p\equiv 1$ mod~$(3)$ and ${\mathcal S}={\mathcal S}_1$ or $\overline{\mathcal S}_1$ then $G<A\cong G_{p,6}\times C_2$. Otherwise, $A=G$.

The surfaces in case~(ix) are the infinite family in case~(v), so they have automorphism groups as described above for case~(v). This completes the proof of Theorem~\ref{integerthm}(c), and thus of Theorem~\ref{integerthm}.


\section{Proof of Theorem~\ref{Jacobians}}\label{JacProof}

Assume first that $\mathcal S$ admits an action of the group $G_{p,r}=\langle a, b\mid a^p=b^r=1, \, bab^{-1}=a^{\omega}\rangle$. Define $m:=(p-1)/r$ and choose $k_1,\ldots, k_m\in{\mathbb Z}_p^*={\mathbb Z}_p\setminus\{0\}$ in such a way that
\[{\mathbb Z}_p^*=\sqcup_{j=1}^m\{k_j, k_j\omega, \ldots, k_j\omega^{r-1}\},\]
where $\sqcup$ denotes disjoint union. Then, by considering the method of {\sl little groups\/} of Wigner and Mackey~\cite[p.~62]{Ser}, we find that $G_{p,r}$ has, up to equivalence, $r$ complex irreducible representations of degree $1$, given by
\[U_l:a\mapsto 1,\quad b\mapsto\xi_r^l\]
for $l=0, \ldots, r-1$, where $\xi_s:=\exp(2\pi i/s)$ for any $s\in\mathbb N$, and $m$ complex irreducible representations of degree $r$, given by
\[V_j: a\mapsto{\rm diag}(\xi_p^{k_j}, \xi_p^{k_j\omega},\ldots, \xi_p^{k_j\omega^{r-1}}), \quad
b\mapsto \left(\,\begin{matrix}0&1&0&\cdots&0\cr 0&0&1&\cdots&0\cr &&&\ddots&\cr 0&0&0&\cdots&1\cr 1&0&0&\cdots&0\cr\end{matrix}\,\right)\]
for $j=1,\ldots, m$.

Choose $r$ mutually distinct integers $t_1,\ldots, t_r\in\{1,\ldots,p-1\}$, and consider
\[P=\langle a\rangle\cong C_p\quad{\rm and}\quad Q_i=\langle a^{t_i}b\rangle\cong C_r\]
for $i=1,\ldots, r$. The dimensions of the vector subspaces of the non-trivial irreducible representations of $G_{p,r}$ fixed by $P$ and $Q_i$ are
\[d_{U_l}^P=d_{V_j}^{Q_i}=1 \quad {\rm and} \quad d_{U_l}^{Q_i}=d_{V_j}^P=0.\]
Thus the collection $\{P, Q_1,\ldots, Q_r\}$ is admissible and therefore, as explained in \S\ref{Jacob}, there exists an abelian subvariety $A$ of $J\mathcal S$ such that
\[J{\mathcal S}\sim J({\mathcal S}/P)\times J({\mathcal S}/Q_1)\times\cdots\times J({\mathcal S}/Q_r)\times A.\]
Each $Q_i$ is conjugate to $Q=\langle b\rangle$ and therefore this isogeny is equivalent to
\[J{\mathcal S}\sim J{\mathcal T}\times (J{\mathcal C})^r\times A\]
where ${\mathcal T}={\mathcal S}/P$ and ${\mathcal C}={\mathcal S}/Q$. The fact that $P$ acts freely on $\mathcal S$ (see Lemma~\ref{Sylowlemma}) implies that $J{\mathcal T}$ is an abelian surface. If the signature of the action of $G$ on $\mathcal S$ is $(m_1,\ldots,m_k)$ then that of the action of $Q$ on $\mathcal S$ is $(\gamma; m_1,\ldots,m_k)$ for some $\gamma\ge 0$. Using the Riemann--Hurwitz formula it is now straightforward to see that $\gamma=(p-1)/r$. Finally, by comparing dimensions one sees that $A=0$, giving the required decomposition.

The remaining cases have already been determined. Indeed, cases (i) and (viii) have been worked out in~\cite[Theorem~1]{IR} and~\cite[Theorem~3]{Rey} respectively, while case~(iv) is a particular subcase of case~(ix) (this is because the central factor $C_2$ here provides no help in obtaining a better decomposition; the problem is that $r=3$ is odd, whereas in cases~(i) and (viii) $r$ is even).


\section{Small primes $p$}\label{smallp}

Although for conciseness we have stated and proved Theorem~\ref{integerthm} only for primes $p\ge 7$, the proof extends, with only minor modifications, to the prime $p=5$. In addition to the cases (vii) and (viii) of Theorem~\ref{integerthm}(a), which are still relevant here, the following exceptional actions arise for integers $\rho\ge 3$:
\begin{itemize}
\item[(a)] one action of a group $G\cong V_{25}\rtimes S_3$ of signature $(2,3,10)$ with $\rho=30$, where $V_{25}:=C_5\times C_5$, restricting to actions of subgroups $V_{25}\rtimes C_3$ of signature $(3,3,5)$ with $\rho=15$, of $V_{25}\rtimes C_2$ of signature $(2,5,10)$ with $\rho=10$, and of $V_{25}$ of signature $(5,5,5)$ with $\rho=5$;
\item[(b)] one action of $G\cong S_5$ of signature $(2,4,6)$ with $\rho=24$, restricting to the action of $G_{5,4}\cong AGL_1(5)$ of signature $(2,2,4,4)$ in Theorem~\ref{integerthm}(a)(vii);
\item[(c)] one action of $G\cong C_5\times S_3$ of signature $(2,10,15)$ with $\rho=6$, restricting to an action of a subgroup $C_{15}$ of signature $(5,15,15)$ with $\rho=3$;
\item[(d)] one action of $G\cong C_{20}$ of signature $(4,5,20)$ with $\rho=4$.
\end{itemize}

For the primes $p=2$ and $3$ one can consult the rich literature on group actions of genus $3$ and $4$ in~\cite{Bro, CI, Kim, KuKo, KK}, for example, together with Conder's lists of group actions in~\cite{ConLarge}.


\section{Small values of $\rho$}\label{smallrho}

Although we have restricted attention to integers $\rho\ge 3$, mainly for simplicity of exposition, the cases $\rho=1$ and $2$ are easily dealt with: in the first case $G\cong C_p$ and all elliptic periods $m_i$ in $\sigma$ are equal to $p$, while in the second case $G\cong C_{2p}$ or $D_p$ and each $m_i=2, p$ or $2p$. (For cyclic and dihedral group actions in general, see~\cite{Har} and~\cite{BCGG} respectively.) However, the results are less uniform than for integers $\rho\ge 3$. It is straightforward to determine the possible signatures $\sigma$ in these two cases. When $\rho=1$ they are the following:
\begin{itemize}
\item $(2;-)$ for any $p$,
\item $(1;2^{[4]})$ for $p=2$,
\item $(1;3,3,3)$ for $p=3$,
\item $(2^{[8]})$ for $p=2$,
\item $(3^{[6]})$ for $p=3$,
\item $(5^{[4]})$ for $p=5$.
\end{itemize}
When $\rho=2$ they are:
\begin{itemize}
\item $(1;2,2)$ for any $p$, $G=C_{2p}$ or $D_p$,
\item $(2,5,5,10)$ for $p=5$, $G=C_{10}$,
\item $(2,6,6,6)$ for $p=3$, $G=C_6$,
\item $(3,3,6,6)$ for $p=3$, $G=C_6$,
\item $(4^{[4]})$ for $p=2$, $G=C_4$,
\item $(2,2,2,4,4)$ for $p=2$, $G=C_4$,
\item $(2,2,2,3,6)$ for $p=3$, $G=C_6$,
\item $(2,2,3,3,3)$ for $p=3$, $G=C_6$ or $D_3$.
\end{itemize}

In each case the Teichm\"uller space of groups $\Gamma(\gamma;m_1,\ldots, m_k)$ has dimension $6\gamma+2k-6>0$, so there is an uncountable family of surfaces $\mathcal S$ admitting the action of $G$. It is a routine matter to count the possible kernels $K$ in a specific group $\Gamma$: for instance, if $\rho=1$, so that $G\cong C_p$, one can use the following results.

\begin{lemma}\label{modsum}
If $p$ is prime and $k\ge 1$ then the number $s_k$ of $k$-tuples $(x_1,\ldots, x_k)\in{\mathbb Z}_p^k$ with each $x_i\ne 0$ and $\sum_{i=1}^kx_i=0$ is given by
\[s_k=\frac{p-1}{p}\left((p-1)^{k-1}+(-1)^k\right).\]
\end{lemma}

\noindent{\sl Proof.} This formula can be proved by applying the Inclusion-Exclusion Principle to the set of {\sl all\/} solutions in ${\mathbb Z}_p$ of $\sum x_i=0$ (without the restriction $x_i\ne 0$), and then excluding those with some $x_i=0$. Alternatively, it can be proved by induction on $k$, using the obvious recurrence relation (consider $x_k$)
\[(s_k, t_k) = (s_{k-1}, t_{k-1})\left(\,\begin{matrix}0&1\cr p-1&p-2\cr \end{matrix}\,\right),
\quad (s_1, t_1) = (0, 1),\]
where
\[t_k=\frac{1}{p-1}\left((p-1)^k-m_k\right)=\frac{1}{p}\left((p-1)^k-(-1)^k\right)\]
is the number of solutions $x_i\ne 0$  in ${\mathbb Z}_p$ of $\sum x_i=a$ for a given $a\ne 0$ (clearly independent of $a$). \hfill$\square$

\begin{cor}\label{counting}
Let $\Gamma$ be a cocompact Fuchsian group with signature $(\gamma;m_1,\ldots, m_k)$. Then the number of normal surface subgroups of prime index $p$ in $\Gamma$ is $0$ unless $m_i=p$ for each $i$, in which case the number is
\[\frac{p^{2\gamma}-1}{p-1}\quad \hbox{if}\quad k=0,\]
and
\[p^{2\gamma-1}\left((p-1)^{k-1}+(-1)^k\right)\quad \hbox{if}\quad k\ge 1.\]
\end{cor}

\noindent{\sl Proof.} The number of such subgroups is equal to the number of surface epimorphisms $\Gamma\to G\cong C_p\cong{\mathbb Z}_p$, divided by the number $p-1$ of automorphisms of $G$. Clearly there are no such epimorphisms unless each $m_i=p$. In this case the epimorphisms correspond bijectively to the choices of elements $a_j, b_j\;(j=1,\ldots,\gamma)$ and $x_i\;(i=1,\ldots, k)$ of ${\mathbb Z}_p$ which generate ${\mathbb Z}_p$ (equivalently are not all zero), with each $x_i\ne 0$ and $\sum x_i=0$. If $k=0$ then any choice of the $2\gamma$ elements $a_j, b_j$, except taking all equal to $0$, is allowed, so the required number is $(p^{2\gamma}-1)/(p-1)$. If $k\ge 1$ then any choice of the elements $a_j, b_j$ is allowed, while Lemma~\ref{modsum} gives the number of choices for the elements $x_i$, so multiplying these leads to the required formula. \hfill$\square$

\medskip

It follows from Corollary~\ref{counting} that in the cases listed above for $\rho=1$, the numbers of normal surface subgroups of index $p$ are $p^3+p^2+p+1$, $4, 9, 1, 11$ and $13$. Similar arguments show that for the signatures listed above for $\rho=2$ the numbers are $4(p+1)$ (for each of $G=C_{2p}$ and $D_p$), $3$, $1$, $3$, $3$, $1$, $1$, $1$ (for $G=C_6$) and $4$ (for $G=D_3$). In many cases, some of these kernels will be conjugate in $PSL_2({\mathbb R})$, leading to isomorphic surfaces $\mathcal S$, and in many cases ${\rm Aut}\,{\mathcal S}$ will be larger than $G$. Lloyd has enumerated equivalence classes of surface epimorphisms $\Gamma\to G\cong C_p$ under the action of ${\rm Aut}\,\Gamma\times{\rm Aut}\, G$ in~\cite{Llo}.


\section{Connections with maps and hypermaps}\label{maps}

Many of the groups $G$ we have classified in Theorem~\ref{integerthm}(a), specifically those in cases (i) to (iv), (vi) and (x) to (xii), arise as quotients of triangle groups $\Gamma=\Gamma(l,m,n)$. As such they are also automorphism groups of orientably regular hypermaps $\mathcal H$, equivalently regular dessins d'enfants (see~\cite{JW}), of type $(l,m,n)$; if $l=2$ these are maps of type $\{m,n\}$ (or dually $\{n,m\}$) in the notation of Coxeter and Moser~\cite{CM}, with $m$ and $n$ the common valencies of the faces and vertices. (More generally, if any of the periods $l, m$ or $n$ is $2$, then by renaming the generators of $\Gamma$ one can regard $\mathcal H$ as a map.) By contrast with the situation we have considered for Riemann surfaces in Theorem~\ref{integerthm}(c), where ${\rm Aut}\,{\mathcal S}\cong N(M)/M$ with $N(M)$ the normaliser of $M$ in $PSL_2({\mathbb R})$, here $G$ is always the full orientation-preserving automorphism group ${\rm Aut}\,{\mathcal H}\cong N_{\Gamma}(M)=\Gamma/M$ of $\mathcal H$, rather than a subgroup of it. Provided the genus $g=p+1$ is not too large, these hypermaps and maps can be found in Conder's computer-generated lists of such objects in~\cite{Con}, or (in the case of maps) in Poto\v cnik's census of rotary maps~\cite{Pot}. By restricting Theorem~\ref{integerthm} to the cases where $\Gamma$ is a triangle group, i.e. ignoring  cases~(v), (vii), (viii) and (ix), and noting that no triangle groups $\Gamma$ arise for $\rho=1$ or $2$ (see Section~\ref{smallrho}), we obtain the classification in Theorem~\ref{mapsthm}, where the numbering of cases follows and refers to that in Theorem~\ref{integerthm}(a).

For small $p$ these maps and hypermaps correspond to entries in Conder's lists of chiral maps, chiral hypermaps, regular maps and regular proper hypermaps in~\cite{Con} as follows (see later in this section for comments on the number of objects represented by each entry, and their chirality, duality and triality properties):

\begin{itemize}
\item[{\rm(i)}] the chiral maps of type $\{6,6\}$ in case~(i) correspond to entry C8.1 for $p=7$, C14.1 for $p=13$, C20.1 for $p=19$, etc;
\item[{\rm(ii)}] the chiral maps of type $\{5,10\}$ in case~(ii) correspond to C12.1 and C12.2 for $p=11$, C32.1 and C32.2 for $p=31$, C42.1 and C42.2 for $p=41$, etc;
\item[{\rm(ii)}] the chiral maps of type $\{8,8\}$ in case~(iii) correspond to C18.1 for $p=17$, C42.3 for $p=41$, etc;
\item[{\rm(iv)}] the chiral hypermaps of type $(3,6,6)$ in case~(iv) correspond to CH8.1 and CH8.2 respectively for $p=7$, CH14.1 and CH14.2 for $p=13$ (with automorphism groups $G_{p,6}$ and $G_{p,3}\times C_2$ respectively in both cases), etc;
\item[{\rm(vi)}] the chiral hypermaps of type $(5,5,5)$ in case~(vi) correspond to CH12.3 and CH12.4 for $p=11$, CH32.3 and CH32.4 for $p=31$, etc;
\item[{\rm(x)}] the regular maps of type $\{3,7\}$ in case~(x) correspond to R14.1, R14.2 and R14.3;
\item[{\rm(xi)}] the regular maps of type $\{3,8\}$ in case~(xi) correspond to R8.1 and R8.2; 
\item[{\rm(xii)}] the regular hypermaps of type $(3,3,4)$ in case~(xii) correspond to RPH8.1 and RPH8.2.
\end{itemize}

For example, in case~(iv) one can distinguish between the two groups $G\cong G_{p,6}$ and $G_{p,3}\times C_2$ by the fact that the former contains $p$ involutions while the latter contains one (necessarily central). Thus for $p=7$, where the only chiral hypermaps of genus $8$ and type $\{3,6,6\}$ listed in~\cite{Con} are CH8.1 and CH8.2, the presentations
\[G_1=\langle R, S\mid R^3=(RS^{-1})^2=S^6=S^{-2}R^{-1}S^{-1}R^{-1}S^{-2}R^{-1}S^{-1}=1\rangle,\]
\[G_2=\langle R, S\mid R^3=S^6=S^{-1}R^{-1}S^{-1}R^{-1}S^2R^{-1}= S^{-1}RS^{-1}R^{-1}SR^{-1}S^{-1}=1\rangle\]
given there for the automorphism groups $G_i$ of CH8.$i$ ($i=1,2$) show that $G_1$ has at least two involutions (namely $RS^{-1}$ and $S^3$, distinct since $G$ is not cyclic), so this must be $G_{7,6}$, while $G_2\cong G_{7,3}\times C_2$ with a unique involution $S^3$. (Note that in $G_2$, if we factor out the central subgroup $C_2$ by putting $S^3=1$, the third defining relation implies that $(RS)^3=1$, so we have a quotient of $\Gamma(3,3,3)$, namely $G_{7,3}$ as expected.) Similar arguments apply for $p=13$ and $19$. However, for some primes $p\equiv 1$ mod~$(3)$, such as $p=31, 37, 43$ and $61$, the numbering of the relevant entries in~\cite{Con} means that the corresponding hypermaps have automorphism groups $G_{p,3}\times C_2$ and $G_{p,6}$ in that reversed order.

Note that in case~(i), for each prime $p\equiv 1$ mod~$(3)$ we found four normal surface subgroups $K$ of $\Gamma=\Gamma(2,6,6)$ with $\Gamma/K\cong G\cong G_{p,6}\times C_2$, representing two chiral pairs of Riemann surfaces $\mathcal S$. However, conjugacy of pairs of subgroups $K$ in the normaliser $N(\Gamma)=\Gamma(2,4,6)$ of $\Gamma$ induces isomorphisms between the two chiral pairs, so up to isomorphism we obtained only one chiral pair of surfaces, ${\mathcal S}_1$ and $\overline{\mathcal S}_1$. Nevertheless, these four subgroups $K$ of $\Gamma$ correspond to four mutually non-isomorphic maps of type $\{6,6\}$, each chiral pair being the vertex-face dual of the other. A similar phenomenon occurs in case~(iii).

More generally, entries in Conder's lists~\cite{Con} represent maps or hypermaps up to chirality and duality (and also triality, interchanging hypervertices, hyperedges and hyperfaces, in the case of hypermaps), so each entry can represent up to four or twelve non-isomorphic maps or hypermaps. Thus entries C12.1 and C12.2 in~\cite{Con}, corresponding to case~(ii) with $p=11$, each represent a chiral pair of maps of type $\{5,10\}$ together with the chiral pair of dual maps of type $\{10,5\}$. Similarly, CH12.3 and CH12.4, corresponding to case~(vi) with $p=11$, each represent six non-isomorphic hypermaps: each set of six is an orbit of the group $C_2\times S_3$ generated by the operations of chirality, duality and triality, induced by the normal inclusion of $\Gamma(5,5,5)$ in the extended triangle group of type $(2,3,10)$ with this quotient.

For the infinite families of chiral maps and hypermaps in cases (i) to (iv) and (vi), $G$ is always the full automorphism group. However, the finitely many exceptional examples in cases (x) to (xii) are all regular, with full automorphism group $A$ containing $G$ with index $2$. In these cases the elements of $A\setminus G$ reverse orientation, and correspond to anticonformal automorphisms of the Riemann surface $\mathcal S$, that is, automorphisms of $\mathcal S$ as a Klein surface. We will now determine these groups $A$.

The regular maps R14.1, R14.2 and R14.3 in case~(x) with $G\cong PSL_2(13)$ are distinguished in~\cite{Con} as having Petrie polygons (closed zigzag paths, turning first right and first left at alternate vertices) of lengths $12$, $26$ and $14$ respectively (this length is twice the order of the commutator $[x,y]$, where $(x,y,z)$ is the canonical generating triple of type $(2,3,7)$ for $G$). An equivalent group-theoretic distinction is that $z$ has trace $t=\pm 6, \pm 5, \pm 3$ respectively, belonging to each of the three conjugacy classes of elements of order $7$ in $G$ (see~\cite[Example 5.4]{JW}). Now Singerman~\cite{Sin74} has shown that any orientably regular map with orientation-preserving automorphism group $G\cong PSL_2(q)$ for some prime power $q$ is in fact regular, with full automorphism group $A\cong PSL_2(q)\times C_2$ or $PGL_2(q)$ as two of the three canonical generators of $G$ are inverted by an inner or outer automorphism of $G$; moreover Hall~\cite[Theorem~2.9]{Hal} has shown that in the case of a Hurwitz group $PSL_2(q)$, these two cases correspond to $3-t^2$ being a square or non-square in ${\mathbb F}_q$, where $t$ is defined as above. The maps R14.1, R14.2 and R14.3 have $3-t^2=6, -1$ and $ -6$, with only $-1$ a square mod~$(13)$, so they have automorphism groups $A\cong PGL_2(13)$, $PSL_2(13)\times C_2$ and $PGL_2(13)$ respectively.

The regular maps R8.1 and R8.2 in case~(xi)  have $G\cong PGL_2(7)$. As the automorphism group of a non-abelian simple group $PSL_2(7)$, $G$ is complete, with no outer automorphisms, so these maps both have automorphism group $A\cong PGL_2(7)\times C_2$ (see~\cite[Exercise~7.17 and Theorem~7.4]{Rot}, for example).

To deal with the regular hypermaps RPH8.1 and RPH8.2 in case~(xii) we need to know the canonical generating triples for their orientation-preserving automorphism groups $G\cong PSL_2(7)$. For $i=1,2$ let ${\mathcal H}_i$ be the orientably regular hypermap of type $(3,3,4)$ corresponding to the following generating triple $(x_i, y_i, z_i)$ for $G$, where matrices in $SL_2(7)$ and $GL_2(7)$ are used to represent elements of $PSL_2(7)$ and its automorphism group $PGL_2(7)$, with the usual convention that scalar matrices represent the identity:
\[x_1, y_1, z_1=\left(\,\begin{matrix}3&0\cr 0&5\cr \end{matrix}\,\right),\;\left(\,\begin{matrix}1&2\cr 3&0\cr \end{matrix}\,\right),\;
\left(\,\begin{matrix}0&1\cr 6&3\cr \end{matrix}\,\right),\]
\[x_2, y_2, z_2=\left(\,\begin{matrix}3&5\cr 0&5\cr \end{matrix}\,\right),\;\left(\,\begin{matrix}1&2\cr 3&0\cr \end{matrix}\,\right),\;
\left(\,\begin{matrix}0&1\cr 6&4\cr \end{matrix}\,\right).\]
For each $i$ the generators $x_i, y_i$ of order $3$ are inverted by conjugation by the involution $g_i$ where
\[g_1=\left(\,\begin{matrix}0&1\cr 2&0\cr \end{matrix}\,\right)\in PGL_2(7)\setminus PSL_2(7)
\quad\hbox{and}\quad g_2=\left(\,\begin{matrix}3&1\cr 4&4\cr \end{matrix}\,\right)\in PSL_2(7).\]
It follows that ${\mathcal H}_1$ and ${\mathcal H}_2$ are both regular, with automorphism groups $PGL_2(7)$ and $PSL_2(7)\times C_2$ respectively. By their genus and type they must be RPH8.1 and RPH8.2 in some order. One can check that the second triple satisfies the defining relation $(acab)^4=1$ for RPH8.2 given in~\cite{Con}, with $x_2=cb$, $y_2=ba$ and $z_2=ac$, whereas the first triple does not, so ${\mathcal H}_i$ is RPH8.$i$ for $i=1, 2$, with $A\cong PGL_2(7)$ and $PSL_2(7)\times C_2$ respectively.

Having mentioned the subject of Petrie length, we note that for the chiral maps in cases (i), (ii) and (iii) the commutator $[x,y]$ is always a non-identity element of $P$, so it has order $p$ and therefore the Petrie length of the map is $2p$.

There are three instances in Theorem~\ref{mapsthm}, namely cases~(i) and (iv), (ii) and (vi), and (xi) and (xii), where the same surfaces $\mathcal S$ support orientably regular maps $\mathcal M$ with orientation-preserving automorphism group $G$, and also orientably regular hypermaps $\mathcal H$ with orientation-preserving automorphism group $G_0$ of index $2$ in $G$. In the second and third of these instances, $\mathcal M$ can be obtained from $\mathcal H$ by representing the latter as its Walsh bipartite map~\cite{Wal} on the same surface, with black and white vertices corresponding to the hypervertices and hyperedges of $\mathcal H$, and edges corresponding to their incidences, and then ignoring the colours of the vertices. The same applies in the first instance, except that here we must first use a triality operation to replace $\mathcal H$ with an orientably regular hypermap of type $(6,6,3)$. In each instance, every automorphism of $\mathcal H$ induces an automorphism of $\mathcal M$, whereas $\mathcal M$ has colour-transposing automorphisms which correspond to dualities rather than automorphisms of $\mathcal H$, thus giving the index~$2$ inclusion between their automorphism groups. All instances of this phenomenon can be explained by the index $2$ inclusion of the triangle group $\Gamma(m,m,n)$ in $\Gamma(2,m,2n)$, see~\cite{Sin72}.

Although Theorem~\ref{mapsthm} is restricted to primes $p\ge 7$, the comments in Section~\ref{smallp} concerning $p=5$ also apply here. The triangle group actions listed there correspond to the following groups $G$, signatures $\sigma$ and entries in~\cite{Con}:
\begin{itemize}
\item[(a)] $V_{25}\rtimes S_3$, $(2,3,10)$, R6.1; $V_{25}\rtimes C_3$, $(3,3,5)$, RPH6.1; $V_{25}\rtimes C_2$, $(2,5,10)$, R6.6; $V_{25}$, $(5,5,5)$, RPH6.11;
\item[(b)] $S_5$, $(2,4,6)$, R6.2;
\item[(c)] $C_5\times S_3$, $(2,10,15)$, R6.10; $C_{15}$, $(5,15,15)$, RPH6.12 and RPH6.13;
\item[(d)] $C_{20}$, $(4,5,20)$, RPH6.7.
\end{itemize}
For example in (b), R6.2 corresponds to the generating triple $((1,2)(3,5), (2,3,4,5), (1,2,3)(4,5))$ of $G=S_5$; the first two generators are inverted by the involution $(3,5)\in G$, so $A=S_5\times C_2$. In (c), entries RPH6.12 and RPH6.13 in~\cite{Con} refer to three regular hypermaps of type $(5,15,15)$ with $G=C_{15}$ and $A=D_{15}$; RPH6.12, corresponding to the generating triple $(3,11,1)$ of $G={\mathbb Z}_{15}$, is a single hypermap, invariant under the duality interchanging hyperedges and hyperfaces (transposing the generators of $\Gamma$ of order $15$), while RPH6.13 consists of a dual pair, corresponding to the triples $(6,8,1)$ and $(12,2,1)$.

With this extension, the results in this section represent a classification of the orientably regular maps and hypermaps of genus $p+1$ with orientation-preserving automorphism group $G$ of order divisible by the prime $p\ge 5$. Much of this (and more, where $p$ does not divide $|G|$) has already been achieved for maps by Conder, \v Sir\'a\v n and Tucker in~\cite{CST}; here we have widened the context to include hypermaps and to relate these combinatorial structures to their underlying Riemann surfaces.


\section{Non-orientable maps and hypermaps}\label{nonor}

If $\mathcal H$ is a non-orientable regular hypermap of type $(l,m,n)$ then its automorphism group $G$ is a quotient $\Delta/M$ of the extended triangle group $\Delta=\Gamma[l,m,n]$ of that type. Its orientable double cover $\tilde{\mathcal H}$ is the orientable regular map corresponding to the map subgroup $\tilde M=M\cap\Gamma$ of the corresponding triangle group $\Gamma=\Gamma(l,m,n)$ (the even subgroup of $\Delta$), with full automorphism group
\[\Delta/\tilde M=(\Gamma/\tilde M)\times(M/\tilde M)\cong(\Delta/M)\times(\Delta/\Gamma)\cong G\times C_2\]
and orientation-preserving automorphism group $\Gamma/\tilde M\cong G$. If $\mathcal H$ has characteristic $-p$ (so that it has genus $p+2$), then $\tilde{\mathcal H}$ has characteristic $-2p$ and hence has genus $p+1$. In particular, if $p$ is a prime dividing $|G|$ and $p\ge 7$ then $\tilde H$ must be one of the regular hypermaps described in Theorem~\ref{mapsthm}, namely one of the three maps R14.1, R14.2 and R14.3 of type $\{3,7\}$ with $p=13$ in case~(x), or one of the two maps R8.1 and R8.2 of type $\{3,8\}$ or hypermaps RPH8.1 and RPH8.2 of type $\{3,3,4\}$ with $p=7$ in cases~(xi) and (xii). We will deal with these possibilities in turn.

As shown in Section~\ref{maps}, the maps R14.1, R14.2 and R14.3 have full automorphism groups $A\cong PGL_2(13)$, $PSL_2(13)\times C_2$ and $PGL_2(13)$ respectively. Only the second of these has the form $G\times C_2$, so we obtain a non-orientable regular quotient map R14.2$/C_2$ of type $\{3,7\}$ and genus $15$ with automorphism group $G\cong PSL_2(13)$. This must be N15.1, the only non-orientable regular map of this type and genus listed in~\cite{Con}. (Note also that  N15.1 has Petrie length $13$; since Petrie lengths are either preserved or halved by factoring out a central subgroup $C_2$, this confirms that N15.1 is not a quotient of R14.1 or R14.3.)

We have seen that R8.1 and R8.2 have full automorphism group $A\cong PGL_2(7)\times C_2$. They therefore yield non-orientable regular quotient maps of type $\{3,8\}$ and genus $9$ with automorphism group $G\cong PGL_2(7)$. These must be N9.1 and N9.2 in~\cite{Con}, in some order. Since these maps have Petrie lengths $7$ and $8$, while R8.1 and R8.2 have Petrie lengths $8$ and $14$, it follows that N9.1 and N9.2 are quotients of R8.2 and R8.1 respectively.

We have also seen that RPH8.1 and RPH8.2 have automorphism groups $A\cong PGL_2(7)$ and $PSL_2(7)\times C_2$ respectively.
We therefore obtain a non-orientable regular hypermap ${\rm RPH}8.2/C_2$ of type $(3,3,4)$ and genus $9$ with automorphism group $PSL_2(7)$. This must be NPH9.1, the only hypermap in~\cite{Con} satisfying this description. However, we obtain no non-orientable hypermap from RPH8.1. Thus we have proved the following (again with the numbering as in Theorem~\ref{integerthm}(a)):

\begin{theo}
The non-orientable hypermaps of characteristic $-p$ for some prime $p\ge 7$ dividing the order of their automorphism group $G$ are (up to triality) as follows:
\begin{itemize}
\item[{\rm (x)}] the regular map {\rm N15.1} of type $\{3,7\}$ and genus $15$ for $p=13$, with orientable double cover {\rm R14.2} in Theorem~\ref{mapsthm}(x), and with $G\cong PSL_2(13)$;
\item[{\rm (xi)}] the regular maps {\rm N9.1} and {\rm N9.2} of type $\{3,8\}$ and genus $9$ for $p=7$, with orientable double covers {\rm R8.2} and {\rm R8.1} in Theorem~\ref{mapsthm}(xi), and with $G\cong PGL_2(7)$;
\item[{\rm (xii)}]  the regular hypermap {\rm NPH9.1} of type $(3,3,4)$ and genus $9$ for $p=7$, with orientable double cover {\rm RPH8.2} in Theorem~\ref{mapsthm}(xii), and with $G\cong PSL_2(7)$.
\end{itemize}
\end{theo}

As in Section~\ref{maps}, it is straightforward to extend this classification to the case $p=5$. The only new example arising is the non-orientable regular map N7.1 of type $\{4,6\}$ and genus $7$, with automorphism group $S_5$ and orientable double cover R6.2 in case~(b) of Section~\ref{maps}. As in the orientable case, these results overlap those in~\cite{CST}, where non-orientable maps are classified without the restriction that $p$ divides $|G|$.


\section{Closing remarks}

We close this paper with an observation and two questions. It is noticeable that the infinite families of maps and hypermaps, corresponding to cases (i) to (iv) and (vi) of Theorem~\ref{integerthm}(a), all occur in chiral pairs, whereas the finitely many sporadic examples, corresponding to cases (x) to (xii) where $p=7$ or $13$, are all (fully) regular, possessing orientation-reversing automorphisms. There is a similar distinction between the underlying Riemann surfaces $\mathcal S$, forming conjugate pairs in cases~(i) to (ix) but not (x) to (xii). In a sense this is partly explained by the method of proof of Theorem~\ref{integerthm}, based on the decomposition (for the all but finitely many $p$ satisfying Lemma~\ref{Sylowlemma}) of the module $M$, which mirrors the decomposition of the homology module $H_1({\mathcal T};{\mathbb C})$ in terms of holomorphic and antiholomorphic differentials. However, a truly satisfactory explanation should depend, not on the human choice of a method of proof, but on intrinsic properties of the objects studied. One might argue that this is just another instance of the well-known phenomenon in finite group theory of infinite families exhibiting uniform behaviour, with finitely many relatively small exceptions, but this does not explain why the infinite families should all be chiral, and the exceptions all regular.

So firstly, is there a better explanation of this phenomenon, and secondly, is it an indication of something more general about the balance between regularity and chirality, or is it simply a consequence of the rather restrictive assumptions applied here?


 \end{document}